\documentstyle[12pt]{article}

\newcommand{\mod}{\mathop{\rm mod}\limits}

\newcommand{\argmax}{\mathop{\rm argmax}\limits}

\newcommand{\const}{\mathop{\rm const}\limits}

\newcommand{\vraisup}{\mathop{\rm vraisup}\limits}

\newcommand{\supp}{\mathop{\rm supp}\limits}

\newcommand{\Law}{\mathop{\rm Law}\limits}

\newcommand{\Ent}{\mathop{\rm Ent}\limits}

\newcommand{\cov}{\mathop{\rm cov}\limits}

\newcommand{\sub}{\mathop{\rm sub}\limits}

\begin{document}

\begin{center}

{\bf Monte-Carlo method for multiple parametric integrals \\
calculation and solving of linear integral Fredholm equations \\
 of a second kind, with confidence regions in uniform norm.}\\

\vspace{3mm}

 $ {\bf E.Ostrovsky^a, \ \ L.Sirota^b } $ \\

\vspace{4mm}

$ ^a $ Corresponding Author. Department of Mathematics and computer science, Bar-Ilan University, 84105, Ramat Gan, Israel.\\
\end{center}
E - mail: \ eugostrovsky@list.ru\\
\begin{center}
$ ^b $  Department of Mathematics and computer science. Bar-Ilan University,
84105, Ramat Gan, Israel.\\
\end{center}
E - mail: \ sirota@zahav.net.il\\

\vspace{3mm}

{\bf Abstract}.  In this article we offer some modification of Monte-Carlo
method for multiple parametric integral computation and
solving of a linear integral Fredholm equation of a second kind (well
posed problem). \par
 We prove that the rate of convergence of offered method is optimal under natural
conditions still in the uniform norm, and construct an asymptotical  and
non-asymptotical confidence region, again in the uniform norm.\par

\vspace{3mm}

 \par

\vspace{3mm}


\vspace{3mm}

\section{Introduction. Notations. Problem Statement. Assumptions.} \par

\vspace{3mm}

 We consider a linear integral Fredholm's equation of a second kind

 $$
 y(t) = f(t) + \int_T K(t,s) \ y(s) \ \mu(ds) = f(t) + S[y](t). \eqno(1.1)
 $$

 Here $ (T = \{t\}, \mu), \ s \in T $ be me measurable space with a probabilistic:
 $ \mu(T) = 1 $ non-trivial measure $ \mu, \ S[y](t) $  is a linear integral operator
 (kernel operator)  with the bimeasurable kernel $ K(\cdot, \cdot):$
$$
S[y](t) \stackrel{def}{=} \int_T K(t,s) \ y(s) \ \mu(ds). \eqno(1.2)
$$
 For example, the set $ T $ may be subset of the whole Euclidean space $ R^d $  with
 non-zero Lebesgue measure $ \mu(ds) = ds. $ \par
  The case when the domain $ T $ dependent on the variable $ t $ may be reduce
in general case after some substitution  to the case of equation (1.1); it may be
implemented, e.g., for the Volterra equation \cite{Grigorjeva1}. \par
  The equations of a view (1.1) appears in many physical problems (transfer equation,
potential theory etc.), in the reliability theory (renewal equation), in the
numerical analysis, for instance, for computation of eigenfunctions and eigenvalues
for integral operators etc. \par
We denote as ordinary for arbitrary measurable function $ g:T \to R $

$$
|g|_p = \left[\int_T |g(s)|^p \ \mu(ds) \right]^{1/p}, \ p \in [1,\infty],
$$

$$
g \in L(p) \ \Leftrightarrow |g|_p < \infty,
$$
where  we will write in the case $ p = \infty $

$$
|g|_{\infty} = \vraisup_{t \in T} |g(t)|.
$$
 We define analogously for the random variable $ \xi $

 $$
 |\xi|_p = \left[{\bf E} |\xi|^p  \right]^{1/p}, p \ge 1.
 $$
\vspace{2mm}
Let us define the following important function, presumed to be finite
$ \mu \ - $ almost everywhere:

$$
R(s) = \vraisup_{t \in T} |K(t,s)|  \eqno(1.3)
$$
and we introduce also the so-called {\it natural } distance, more exactly, semi-distance,
 $ d = d(t,s) $ on the space $ T: $

$$
d(t,s) \stackrel{def}{=} \vraisup_{x \in T} \frac{|K(t,x)-K(s,x)|}{R(x)}, \eqno(1.4)
$$
so that

$$
|K(t,x) - K(s,x)| \le R(x) \ d(t,s). \eqno(1.5)
$$

We assume that the metric space $ (T,d) $ is compact set and the measure
$ \mu $  is Borelian.  We suppose also the  function $ f(t) $ is $ d \ - $ continuous
$ f: T \to R^1. $ \par
It follows from the inequality (1.5) that the function $ t \to K(t,x) $ is also
$ d \ - $ continuous  for $ \mu \ - $ almost everywhere values $ x; \ x \in T. $ \par
 For the integral operator $ S[\cdot] $  the $ m^{th}, m \ge 2 $ power of
 $ S: S^m[\cdot] $ may be calculated as usually

 $$
 S^m[g](t) = \int_{T^m} K(t,s_1)\prod_{j=1}^{m-1}K(s_j,s_{j+1}) \ g(s_m) \
 \prod_{i=1}^m \mu(ds_i). \eqno(1.6)
 $$

 We will consider the source equation in the space $ C(T) = C(T,d) $ of all $ d \ - $ continuous numerical functions $ g: T \to R $  with ordinary norm

 $$
 ||g|| = \sup_{t \in T} |g(t)| = |g|_{\infty}. \eqno(1.7).
 $$

 Recall that the norm of linear integral operator $ S[\cdot] $  (1.2) in this space may
 be calculated by the formula

 $$
 ||S|| \stackrel{def}{=} \sup_{g \in C(T), g \ne 0} ||S[g]||/||g|| =
 \sup_{t \in T} \int_T |K(t,s)| \ \mu(ds). \eqno(1.8)
 $$

 A spectral radius $ r = r(L) $ of a bounded linear operator $ L: C(T) \to C(T) $ is defined by the expression

$$
r(L) = \overline{\lim}_{m \to \infty} ||L^m||^{1/m} =
\lim_{m \to \infty}||L^m||^{1/m}. \eqno(1.9)
$$
 More detail information about spectral radius of operator see in the classical books
 of N.Dunford, J.Schwartz \cite{Dunford1}, chapter VII, section 3;
 \cite{Dunford2}, chapter IX, section 1.8. For instance, it is known that $ r(L) \le ||L||. $\par
 We define for the kernel operator $ S[\cdot] $ of a view (1.2) the so-called Kroneker's
power $ U = S^{(2)} $ as an kernel integral operator by the following way:

$$
U[g](t) = \int_T K^2(t,s) \ g(s) \ \mu(ds). \eqno(1.10)
$$

We assume (the essential condition!) that

$$
\rho \stackrel{def}{=} r(U) = r \left(S^{(2)} \right) < 1. \eqno(1.11)
$$
 Note that in the considered case (1.11)

 $$
 \rho_1 \stackrel{def}{=} r(S) \le \sqrt{\rho} < 1. \eqno(1.12)
 $$

{\bf Example 1.}  Let the set $ T $ be closed interval $ T = [0,1] $ equipped with
classical Lebesgue measure. Suppose the function $ (t,s) \to K(t,s) $ is continuous and denote

$$
\gamma = \gamma(K) = \sup_{t,s \in T} |K(t,s)|. \eqno(1.13)
$$
If $ \gamma(K) < 1,$  then evidently $ \rho_1(S) \le \gamma < 1, \rho(S) \le
\gamma^2 < 1.  $ \par
{\bf Example 2.} Let us consider instead operator $ S $ the Volterra's operator of a
view

$$
V[g](t) = \int_0^t K(t,s) \ g(s) \ ds,
$$
where $ T = [0,1]$ and we suppose that the function $ (t,s) \to K(t,s) $ is continuous and we will use again the notations (1.13). Then

$$
||V^m|| \le \frac{\gamma^m}{m!}.
$$
Therefore

$$
\rho_1(V) = \rho(V) = 0.
$$
The (linear) operators $ V $ with the property $ r(V) = 0 $ are called quasinilpotent. \par

 Notice that the Monte-Carlo method with optimal rate of convergence for
 linear Volterra's equation

 $$
 y(t) = f(t) + \int_0^t K(t,s) \ y(s) \ ds
 $$
is described, with an applications in the reliability theory, in an article \cite{Grigorjeva1}. We represent in this article  some generalization  of results
of \cite{Grigorjeva1}; but we note that in this article is considered also the case
of {\it discontinuous} function $ f(\cdot, \cdot) $  and are described applications
in the reliability theory (periodical checking, prophylaxis). \par
{\bf Example 3.} We denote for arbitrary linear operator $ L $

$$
 r_m(L) \stackrel{def}{=} ||L^m||. \eqno(1.14)
$$
 If

 $$
 \exists \beta_1 \in (0,1), C_1, \Delta_1 = \const < \infty \ \Rightarrow
 r_m(S) \le C_1 m^{\Delta_1} \beta^m_1, \ m = 1,2,\ldots, \eqno(1.15)
 $$
then obviously  $ \rho_1 \le \beta_1. $ Analogously, if

 $$
 \exists \beta \in (0,1), C, \Delta = \const < \infty \ \Rightarrow
 r_m(U) \le C m^{\Delta} \beta^m, \ m = 1,2,\ldots, \eqno(1.16)
 $$
then  $ \rho \le \beta. $ \par

\vspace{3mm}
{\bf We intend to prove that under formulated conditions that there exists a Monte-Carlo
method for solving of equation (1.2) with the classical speed of convergence
$ 1/\sqrt{n}, $  where $ n $ denotes the common number of elapsed random variables.\par
Moreover, at the same result is true when the convergence is understudied in the
uniform norm. }\par

\vspace{3mm}

The letter $ C, $ with or without subscript, denotes a finite positive non essential constants,  not necessarily  the same at each appearance. \par

\vspace{3mm}

 The papier is organized as follows. In the next section we describe the numerical
Monte-Carlo  method for solving of integral equation (1.1) and prove the optimality
of it convergence in each fixed point $ t_0; \ t_0 \in T. $ \par
In the third section we recall for reader convenience
some used facts about Grand Lebesgue Spaces of random
variables and random processes and obtain some new used results. \par
 Fourth section is devoted to the Monte-Carlo computation of multiple parametric
integrals.  The fifth section contains the main result of offered article: confidence
region for solution of Fredholm's integral equation in the uniform norm under the
classical normalizing.
In the next section we built the non-asymptotical confidence domain for
multiple parametric  integrals  and for solution of Fredholm's integral equation
in the uniform norm. In the $ 7^{th} $ section we consider some examples to show
the convenience of using of offered algorithms. \par
 In the next section we offer the Monte-Carlo method for derivative computation
for solution of Fredholm's integral equation again in the uniform norm with optimal
rate of convergence.\par
 The $ 9^{th} $ section included some additional remarks.
The last section contains some results about necessity of conditions of our theorems.\par

\vspace{3mm}

\section{Numerical method.  Speed of convergence.}

\vspace{3mm}
{\bf  A. Deterministic part of an error.} \par
Let $ \epsilon > 0 $ be arbitrary "small" positive number. The solution of the equation
(1.1) may be written by means of the Newman's series

$$
y(t) = f(t) + \sum_{m=1}^{\infty} S^m[f](t).\eqno(2.1)
$$

Let $ \epsilon \in (0, 1/2) $ be a fix "small" number. We introduce as an approximation
for the solution $ y(\cdot) $ a finite sum

$$
y^{(N)}(t) = \sum_{m=1}^{N} S^m[f](t), \eqno(2.2)
$$
where the amount of summands $ N = N(\epsilon) $ may be determined from the condition

$$
\sum_{m=N+1}^{\infty} ||S^m[f]|| < \epsilon. \eqno(2.3)
$$
If for example the operator $ S[\cdot] $ satisfies the condition (1.15), then the number
$  N = N(\epsilon) $ may be find as the minimal integer solution of the inequality

$$
C_1||f||\sum_{N+1}^{\infty} m^{\Delta_1} \beta^m_1 \le \epsilon \eqno(2.4)
$$
in the domain
$$
N \ge \argmax_m \left[ m^{\Delta_1} \beta^m_1 \right] = \frac{\Delta_1}{|\log \beta_1|}.
$$
 The asymptotical as $  \epsilon \to 0+ $ solution $ N \sim N_0 = N_0(\epsilon) $
 of the equation (2.4) has a view

 $$
 N_0 = \log (C_1|\log \beta_1|/\epsilon_1) +
 \Delta_1 \left[\log (C_1|\log \beta_1|/\epsilon_1)/|\log \beta_1| \right],\eqno(2.5)
 $$

$$
\epsilon_1 = \epsilon/\left(C_1 \ ||f|| \right).
$$

\vspace{3mm}
{\bf  B. Probabilistic part of an error.} \par
\vspace{3mm}
We offer for the calculation of the value $ y^{(N)}(t) $ at the {\it fixed}
point $ t \in T $ the Monte-Carlo method. Let us introduce a following notations.\par

$$
R_{\alpha} = R_{\alpha}(U) = \sum_{k=1}^{\infty} k^{\alpha} r_k^{1/2}(U), \
\alpha = \const \in (-\infty,\infty);
$$

$$
 R_{\alpha}(N,U) = \sum_{k=1}^{N} k^{\alpha} r_k^{1/2}(U), \
\alpha = \const \in (-\infty,\infty).
$$
 Further, for the $ m - $ tuple $ \vec{x} = (x_1,x_2, \ldots, x_m), \ m = 2,3, \ldots $ we define

$$
\vec{K}^{(m)}[f](t,\vec{x})= K(t,x_1)K(x_1,x_2)\ldots K(x_{m-1},x_m) f(x_m).\eqno(2.6)
$$
We define in the case $ m=1 $

$$
\vec{K}^{(1)}[f](t,\vec{x})= K(t,x_1)f(x_1).\eqno(2.7)
$$

Let $ \xi_i^{(j)} $ be a double sequence independent $ \mu $ distributed random
variables:

$$
{\bf P} (\xi_i^{(j)} \in A) = \mu(A),
$$
where $ A $ is arbitrary Borelian subset of the space $ T. $ We introduce the following
random vector for the integer values $ m = 1,,2,\ldots: $

$$
\vec{\xi}^{(j)}_m  = \{\xi_1^{j},\xi_2^{j}, \ldots, \xi_m^{j}\}.
$$
 Let us denote
 $$
 \hat{n} = (n(1), n(2), \ldots, n(N)), \
 B(\hat{n}) \stackrel{def}{=} \sum_{m=1}^N m \cdot n(m),
 $$
where $ N  = N(\epsilon), \ \hat{n} $
be any $ N - $ tuple of positive integer numbers:
$ n(j) = 1,2,\ldots. $  We consider the following Monte-Carlo approximation for the
multiple integral $ S^m[f]: $

$$
S^m_{n(j)}[f]  = S^m_{n(j)}[f](t) =
 \frac{1}{n(j)} \sum_{l=1}^{n(j)} \vec{K}^{(m)}[f](t,\vec{\xi}^{(j)}_l)\eqno(2.8)
$$
and  we offer correspondingly the following approximation  $ y^{(N)}_{\hat{n}}(t) $
 for the variable $ y^{(N)}(t): $
$$
y^{(N)}_{\hat{n}}(t):= f(t) + \sum_{j=1}^N S^m_{n(j)}[f](t). \eqno(2.9)
$$
 Note that the approximation (2.8) for the (multiple) parametric integrals was
 introduced by  Frolov A.S.and Tchentzov N.N., see \cite{Frolov1}, and was named
 "Dependent Trial Method".  It was proved under some
 hard conditions that the rate of convergence of this approximation in the space of
 continuous functions is optimal, i.e. coincides with the expression $ 1/\sqrt{n}, $
where $ n $ denotes  the amount of all used random variables. \par
{\bf Theorem 1.1} We assert under formulated conditions that for arbitrary
$ \epsilon \in (0,1) $ and for all sufficiently great values $ n; \ n \to \infty $
there exists the tuple $ \hat{n} = \hat{n}(n) $ for which

$$
\sup_{t \in T} \min_{\hat{n}: B(\hat{n}) \le n}
{\bf Var} \left(y^{(N(\epsilon))}_{\hat{n}}(t) - y^{(N)}(t) \right) \le
$$

$$
||f||^2 \  R_{1/2}(U) \ R_{-1/2}(U) \ \left[ \frac{1}{n} + \frac{C_1(U)}{n^2}  \right], \eqno(2.10)
$$
and

$$
\sup_{t \in T} \min_{\hat{n}: B(\hat{n}) \le n}
{\bf Var} \left(y^{(N(\epsilon))}_{\hat{n}}(t) - y^{(N)}(t) \right) \ge
$$
$$
||f||^2 \  R_{1/2}(U) \ R_{-1/2}(U) \ \left[ \frac{1}{n} - \frac{C_2(U)}{n^2}  \right],
$$
where the positive finite constants
$ C_1= C_1(U), C_2 = C_2(U) $ does not depend on the $ n,\epsilon $ and $ f. $ \par
{\bf  Remark 2.1.} Note that under condition $ \rho(U) $ the values
$  R_{1/2}(U) $ and $ R_{-1/2}(U)  $ are finite. \par
{\bf Proof} of the theorem 1. Without loss of generality we can suppose $ ||f|| = 1.$ \par
We conclude that the variance of each summand in (2.9), i.e. the expression

$$
v(m, n(j))(t) := {\bf Var} \left( S^m_{n(j)}[f](t) \right)
$$
may be estimated as follows:

$$
\sup_{t \in T} v(m, n(j))(t) \le r_m /n(m), \ m = 1,2,\ldots.
$$
Therefore,

$$
\sup_{t \in T} {\bf Var} \left(y^{(N(\epsilon))}_{\hat{n}}(t) \right) \le
\sum_{m=1}^N \frac{r_m}{n(m)} =: \Phi(\hat{n}). \eqno(2.11)
$$

 On the other hand, the common amount of used random variables $ \xi_i^j $ in the
 formula (2.9) is equal to the expression

 $$
 B(\hat{n}) = \sum_{m=1}^N  m \cdot n(m). \eqno(2.12)
 $$
 Let us consider the following constrained extremal problem:

$$
\Phi(\hat{n}) \to \min / B(\hat{n}) = n. \eqno(2.13)
$$
 We obtain using the Lagrange's factors method the optimal value
 $ \hat{n_0}  =  \{n_0(1), n_0(2), \ldots,   \}$  of the tuple $ \hat{n} $ for the
 problem (2.13) has a view:

 $$
 n_0(m) = \frac{n r^{1/2}_m(U)}{R_{1/2}(N,U) \ \sqrt{m}} \stackrel{def}{=}
 \theta(m) \ n,
 $$
up to rounding to the integer number, for example,

$$
n_0(m) =  1+ \Ent
\left[ \frac{n r^{1/2}_m(U)}{R_{1/2}(N,U) \ \sqrt{m}} \right] = \eqno(2.14)
$$
$$
1 + \Ent \left[\theta(m) \ n \right],
$$

where $ \Ent(a) $ denotes the integer part of the positive number $ a $  and

$$
\theta(m)= \theta(m,N,U) = \frac{ r^{1/2}_m(U)}{R_{1/2}(N,U) \ \sqrt{m}}
$$

The minimal  value of the functional
 $ \Phi(\hat{n}) $ under condition $  B(\hat{n}) = n $ may be estimated as

 $$
 \min \Phi(\hat{n})/[B(\hat{n}) = n] \ \le
 R_{1/2}(N,U) \ R_{-1/2}(N,U) \ \left[ \frac{1}{n} + \frac{C_1(U)}{n^2}  \right] \le
 $$

$$
R_{1/2}(U) \ R_{-1/2}(U) \ \left[ \frac{1}{n} + \frac{C_1(U)}{n^2}  \right].
$$
The lower bound provided analogously.\par
 This completes the proof of theorem 1.1. \par
{\bf Example 2.1.} Assume that the sequence $ r_k(S) $ satisfies the condition (1.15).
 As long as as $ x \to 1 - 0 $ and $  \beta = \const > -1 $

$$
\sum_{k=1}^{\infty} k^{\beta} x^{k/2} \sim
\frac{2^{\beta+1} \Gamma(\beta+1)}{|\log x|^{\beta+1}},
$$
we have as $ \alpha = \const > -1/2, \beta_1 \to 1-0 $

$$
\sup_{t \in T} {\bf Var} \left(y^{(N(\epsilon))}_{\hat{n}}(t) - y^{(N)}(t) \right) \sim
$$

$$
C_1 ||f||^{2} \frac{2^{2\alpha+2} \Gamma(\alpha+3/2) \Gamma(\alpha+1/2)}
{|\log \beta_1|^{2\alpha+2}}
 \ \left[ \frac{1}{n} + \frac{C_1(U)}{n^2}  \right].
$$

\vspace{3mm}

{\bf Remark 2.1.} Note that  in the case $ T \subset R^d $ for each random number
$ \vec{\xi} $ generation are used in general case $  d $ uniform distributed in
the interval $ [0,1] $ random variables. See, e.g. \cite{Devroye1}, \cite{Grigorjeva1}.\par

\section{Banach spaces  of random variables}

\vspace{2mm}
{\bf Pilcrow A. } {\it  Banach spaces of random variables with exponentially decreasing
tails of distributions.} ("Exponential" level). \par
 In order to formulate our results, we need to introduce some addition
notations and conditions. Let $ \phi = \phi(\lambda), \lambda \in (-\lambda_0, \lambda_0), \ \lambda_0 = \const \in (0, \infty] $ be some even strong convex which takes positive values for positive arguments twice continuous
differentiable function, such that
$$
 \phi(0) = 0, \ \phi^{//}(0) > 0, \ \lim_{\lambda \to \lambda_0} \phi(\lambda)/\lambda = \infty.
$$
 We denote the set of all these function as $ \Phi; \ \Phi =
\{ \phi(\cdot) \}. $ \par
 We  say that the {\it centered} random variable (r.v) $ \xi = \xi(\omega) $
belongs to the space $ B(\phi), $ if there exists some non-negative constant
$ \tau \ge 0 $ such that

$$
\forall \lambda \in (-\lambda_0, \lambda_0) \ \Rightarrow
{\bf E}\exp(\lambda \xi) \le \exp[ \phi(\lambda \ \tau) ]. \eqno(3.1).
$$
 The minimal value $ \tau $ satisfying (4) is called a $ B(\phi) \ $ norm
of the variable $ \xi, $ write
 $$
 ||\xi||B(\phi) = \inf \{ \tau, \ \tau > 0: \ \forall \lambda \ \Rightarrow
 {\bf E}\exp(\lambda \xi) \le \exp(\phi(\lambda \ \tau)) \}. \eqno(3.2)
 $$
  For instance, if $ \phi(\lambda) \stackrel{def}{=}\phi_2(\lambda) = 0.5 \lambda^2,
  \ \lambda \in R $ the space $ B(\phi_2) $ is called  subgaussian space and
is denoting ordinary $ B(\phi_2) = \sub = \sub(\Omega)  $  in accordance
  to Kahane \cite{Kahane1}; the (centered) random variables from this space are
  called subgaussian.\par
   The norm in subgaussian space  of a random variable $ \xi $
 will denoted $ ||\xi||\sub.  $ \par
 The important example of subgaussian random variables (r.v.) are centered Gaussian
(normal) variables; indeed, if r.v. $ \Law(\xi) = N(0,\sigma^2), \sigma \ge 0, $ then
$ ||\xi||\sub = \sigma. $ \par
 If a centered r.v. $ \xi $ is bounded, then it is also subgaussian and
 $$
 ||\xi||\sub \le  |\xi|_{\infty} := \vraisup |\xi|.
 $$
  For instance, the Rademacher's r.v. $ \xi: $
  $$
  {\bf P}(\xi=1) = {\bf P}(\xi=-1) = 1/2
  $$
is also subgaussian and $ ||\xi||\sub = 1. $ \par
 It is proved in the article \cite{Buldygin2} that the space $ \sub(\Omega) $
is Banach space. The centered random variable $ \xi $ belongs  to the space
$ \sub(\Omega) $ and has a norm $ \tau = ||\xi||\sub, \ \tau \ge 0 $ if and only if

$$
\forall \lambda \in R \ \Rightarrow {\bf E} \exp(\lambda \xi) \le
\exp \left(0.5 \lambda^2 \tau^2 \right).
$$
 More details about the space $ \sub(\Omega) $ see in the article \cite{Buldygin2}.\par
\vspace{3mm}

 The spaces $ B(\phi) $ are rearrangement invariant in the terminology of a book \cite{Bennet1},  chapter 2 and 3; \cite{Krein1}, chapter 3;
  are very convenient for the investigation of the r.v. having a
exponential decreasing tail of distribution, for instance, for investigation of
the limit theorem, the exponential bounds of distribution for sums of random variables,
theory of martingales, non-parametrical statistics,
non-asymptotical properties, problem of continuous of random fields,
study of Central Limit Theorem in the Banach space etc., see \cite{Buldygin1},
\cite{Kozachenko1}, \cite{Talagrand1}, \cite{Talagrand2}, \cite{Talagrand3},
\cite{Talagrand4}, \cite{Ostrovsky1}, \cite{Ostrovsky2},
 \cite{Sirota1}, \cite{Sirota2},\cite{Sirota3},\cite{Sirota4},\cite{Sirota5},
 \cite{Sirota6},\cite{Sirota7},\cite{Sirota8}, \cite{Talagrand1}. \par

 The generalization of this spaces on the case $ \mu(X) = \infty $ is considered,
 e.g., in  \cite{Fiorenza1}, \cite{Fiorenza2}, \cite{Fiorenza3},
  \cite{Iwaniec1}, \cite{Iwaniec2}, \cite{Liflyand1}, \cite{Ostrovsky16}, \cite{Ostrovsky17} etc.

  The space $ B(\phi) $ with respect to the norm $ || \cdot ||B(\phi) $ and
ordinary operations is a Banach space which is isomorphic to the subspace
consisted on all the centered variables of Orlich's space $ (\Omega,F,{\bf P}), N(\cdot) $ with $ N \ - $ function

$$
N(u) = \exp(\phi^*(u)) - 1,
$$
where
$$
 \phi^*(u) \stackrel{def}{=} \sup_{\lambda} (\lambda u -\phi(\lambda)).
$$
 The transform $ \phi \to \phi^* $ is called Young-Fenchel transform. The proof of considered assertion used the properties of saddle-point method and theorem
 of Fenchel-Moraux:
$$
\phi^{**} = \phi,
$$
see  \cite{Krasnoselsky1}, chapter 1. \par
 Let $ \xi $ be centered random variable such that its {\it moment generating function}

 $$
  \lambda \to {\bf E} \exp(\lambda \xi)
 $$
is finite in some neighborhood of origin: $ |\lambda| < \lambda_0, \ \lambda_0 = \const,  0 < \lambda_0 \le \infty. $ \par
 Here $ \lambda $ may be complex; in this case the moment generating function is
 analytical inside the circle $ |\lambda| < \lambda_0. $ \par
 The finiteness of moment generating function is equivalent the following moment
 inequality:

 $$
|\xi|_p \le C \ p, \ p \ge 1.
 $$

The {\it natural  function} $ \phi = \phi_{\xi}(\lambda) $ for the variable $ \xi $
may be introduced by a formula

$$
\phi_{\xi}(\lambda) = \log {\bf E} \exp(\lambda \xi).
$$
 It is obvious that $ \phi_{\xi}(\cdot) \in \Phi. $\par

\vspace{3mm}

 Analogously is defined a so-called co-transform $ v \to v_*: $

$$
v_*(x) \stackrel{def}{=} \inf_{y \in (0,1)} (xy + v(y)).
$$

 The next facts about the $ B(\phi) $ spaces are proved, for instance,
  \cite{Kozachenko1}, \cite{Ostrovsky1},  p. 19-40:

$$
{\bf 1.} \ \xi \in B(\phi) \Leftrightarrow {\bf E } \xi = 0, \ {\bf and} \ \exists C = \const > 0,
$$

$$
U(\xi,x) \le \exp(-\phi^*(Cx)), x \ge 0,
$$
where $ U(\xi,x)$ denotes in this article the {\it tail} of
distribution of the r.v. $ \xi: $

$$
U(\xi,x) = \max \left( {\bf P}(\xi > x), \ {\bf P}(\xi < - x) \right),
\ x \ge 0,
$$
and this estimation is in general case asymptotically exact. \par
 Here and further $ C, C_j, C(i) $ will denote the non-essentially positive
finite "constructive" constants.\par
 More exactly, if $ \lambda_0 = \infty, $ then the following implication holds:

$$
\lim_{\lambda \to \infty} \phi^{-1}(\log {\bf E} \exp(\lambda \xi))/\lambda =
K \in (0, \infty)
$$
if and only if

$$
\lim_{x \to \infty} (\phi^*)^{-1}( |\log U(\xi,x)| )/x = 1/K.
$$
 Here and further $ f^{-1}(\cdot) $ denotes the inverse function to the
function $ f $ on the left-side half-line $ (C, \infty). $ \par
 Let $ \xi = \xi(t) = \xi(t,\omega) $ be some centered random field, $ t \in T. $
 The function $ \phi(\lambda) = \phi_{\xi}(\lambda) $ may be "constructive" introduced
 by the formula
$$
\phi(\lambda) = \phi_0(\lambda) \stackrel{def}{=} \log \sup_{t \in T}
 {\bf E} \exp(\lambda \xi(t)), \eqno(3.3)
$$
 if obviously the family of the centered r.v. $ \{ \xi(t), \ t \in T \} $ satisfies the {\it uniform } Kramer's condition:
$$
\exists \mu \in (0, \infty), \ \sup_{t \in T} U(\xi(t), \ x) \le \exp(-\mu \ x),
\ x \ge 0.
$$
 In this case, i.e. in the case the choice the function $ \phi(\cdot) $ by the
formula (3.3), we will call the function $ \phi(\lambda) = \phi_0(\lambda) $
a {\it natural } function for the random process $ \xi(t).$  \par
\vspace{2mm}
 {\bf Pilcrow B.} {\it Grand Lebesgue Spaces of random Variables. } ("Power" level.) \par
 \vspace{2mm}

 Let $ \psi = \psi(p) $ be function defined on some {\it semi-open} interval of a view

$$
1 \le p < b,
$$
where obviously  $ \ b = \const,  1 < b \le \infty, $ is  continuous and is bounded from
below: $ \inf_{p \in (1,b)} \psi(p) > 0 $ function such that the function

$$
w(p) = w_{\psi}(p) = p \log \psi(p), p \in (1,b)
$$
is downward convex.  We will denote the set of all such a functions  as $ \Psi =
\Psi(1,b): \ \psi \in \Psi. $ \par
 We define as ordinary

 $$
 \supp \psi = [1,b)
 $$
and put for the values $ p \ge b \ \rightarrow \psi(p) = + \infty $ in the case
$ b < \infty. $ \par
 It may be considered analogously the case of closed interval
 $ p \in [1,b], \ b \in (1,\infty];$  then $ \supp \psi = [1,b]  $ and for all the
 values $ p > b \ \rightarrow \psi(p) = + \infty. $  \par
 We introduce a new norm (the so-called "moment norm")
on the set of r.v. defined in our probability space by the following way: the space
$ G(\psi) $ consist, by definition, on all the centered r.v. with finite norm

$$
||\xi||G(\psi) \stackrel{def}{=} \sup_{p \in \supp \psi} |\xi|_p/\psi(p), \ |\xi|_p
\stackrel{def}{=}{\bf E}^{1/p} |\xi|^p. \eqno(3.4)
$$
{\bf Remark 3.1.} Note that in the case $ \supp \psi = [1,b],
\ b = \const \in (1,\infty) $
the norm $ || \xi ||G(\psi) $ coincides with the ordinary Lebesgue norm $ |\xi|_b: $

$$
||\xi||G(\psi) = |\xi|_b.
$$
 Indeed, the inequality $|\xi|_b \le ||\xi||G(\psi)  $ is evident; the inverse inequality
 follows from the Lyapunov's inequality.\par
  Let us introduce the function
$$
\chi(p) = \chi_{\phi}(p) = \frac{p}{\phi^{-1}(p)}, \ p \ge 1.
$$

 It is proved, e.g. in \cite{Kozachenko1},
 \cite{Ostrovsky1}, chapter 1, section (1.8) that
 the spaces $ B(\phi) $ and $ G(\chi) $ coincides:$ B(\phi) =
G(\chi) $ (set equality) and both
the norm $ ||\cdot||B(\phi) $ and $ ||\cdot||G(\chi) $ are equivalent: $ \exists C_1 =
C_1(\phi), C_2 = C_2(\phi) = \const \in (0,\infty), \ \forall \xi \in B(\phi) $

$$
||\xi||G(\chi) \le C_1 \ ||\xi||B(\phi) \le C_2 \ ||\xi||G(\chi).
$$
 Conversely, for arbitrary function $ \psi(\cdot) \in \Psi(1,\infty) $ for which

 $$
 \overline{\lim}_{p \to \infty} \frac{\log \psi(p)}{\log p} < 1
 $$
may be defined the correspondent function $ \phi = \phi(\lambda) $ as follows:
$$
\phi_{\psi}(\lambda) = \left[ \frac{\lambda}{\psi(\lambda)} \right]^{-1}, \
\lambda \ge \lambda_0 = \const > 0.
$$
 Recall that at $ |\lambda| \le \lambda_0 \ \Rightarrow \phi_{\psi}(\lambda) \asymp
 C \lambda^2. $\par
 The definition (3.4) is correct still for the non-centered random
variables $ \xi.$ If for some non-zero r.v. $ \xi \ $ we have $ ||\xi||G(\psi) < \infty, $ then for all positive values $ u $

$$
{\bf P}(|\xi| > u) \le 2 \ \exp \left( - w^*_{\psi}(u/(C_3 \ ||\xi||G(\psi))) \right).
\eqno(3.5)
$$
and conversely if a r.v. $ \xi $ satisfies (3.5), then $ ||\xi||G(\psi) <
\infty. $ \par
The definition (3.4) is more general as (3.1).  Indeed, if a r.v. $ \xi $ belong to some
space $ B(\phi), \ \phi \in \Phi, $ then
$ \forall p \in (1,\infty) \ |\xi|_p < \infty.  $ The inverse inclusion is not true, e.g., for the symmetrical distributed random variable $ \zeta $  for which

$$
{\bf P} (|\zeta| > u) = \exp \left(-u^{\Delta} \right), u \ge 0,  \eqno(3.6)
$$
where $ \Delta = \const \in (0,1). $ \par
 Further, let $ \xi $ be any r.v. such that for some $ b = \const > 1 |\xi|_b < \infty. $
The {\it natural choice} of the function $ \psi_{\xi}(p) $ for the r.v.$ \xi $ may be
defined by the formula

$$
\psi_{\xi}(p) = |\xi|_p, \ p: |\xi|_p < \infty. \eqno(3.7)
$$
{\bf Remark 3.2.} Note that:\par
{\bf A. } The r.v. $ \xi $ is bounded if and only if

$$
\overline{\lim}_{p \to \infty} \psi_{\xi}(p) < \infty.
$$

{\bf B. } The r.v. $ \xi $ satisfies the Kramer's condition if and only if

$$
\overline{\lim}_{p \to \infty} \log \psi_{\xi}(p)/\log p \le 1.
$$

{\bf C. } The r.v. $ \xi $ obeys all the exponential moments, i.e.

$$
\forall \lambda \in R \ \Rightarrow {\bf E} \exp(\lambda \xi) < \infty
$$
if and only if

$$
\overline{\lim}_{p \to \infty} \log \psi_{\xi}(p)/\log p < 1.
$$

\vspace{3mm}

 For instance,  for the r.v. $ \zeta $ in (3.6) the natural function $ \psi_{\zeta}(p) $
has a view

$$
\psi_{\zeta}(p) = |\zeta|_p = 2^{1/p} \Gamma^{1/p} (p/\Delta + 1), \ 1 \le p < \infty.
$$
Note that as $ p \to \infty $

$$
\psi_{\zeta}(p) \sim (p/(\Delta e)^{1/\Delta}.
$$

{\bf Pilcrow C}. {\it Non-asymptotical bounds of distributions in the classical CLT.}\par

Let us define for all the functions $ \phi \in \Phi $

$$
\overline{\phi}(\lambda) = \sup_{n=1,2,\ldots} n \phi(\lambda/\sqrt{n}).
$$
 For example, let

 $$
 \phi(\lambda) = \phi_r(\lambda), \ |\lambda| \ge 1 \Rightarrow
 \phi_r(\lambda) = C_1 |\lambda|^r, \ r = \const > 1;
 $$
then

$$
\overline{\phi_r}(\lambda) \asymp \phi_{\max(r,2)}(\lambda), \ |\lambda| \ge 1.
$$
\vspace{2mm}

 We denote also for all the functions $ \psi \in \Psi(2,b) $

$$
\overline{\psi}(p) =   C_0^{-1}  \cdot  p \ \psi(p)/\log p.
$$
 The provenance, calculation and  exact value of the constant
 $ C_0 \approx 1.77638\ldots  $ is described in \cite{Ostrovsky17}. \par
  Obviously, if $ b < \infty, $ then $ \overline{\psi}(p)\asymp \psi(p). $ It is not true
 in the case when $ b = \infty. $ \par
 The probabilistic sense of introduced functions is following. Let $ \xi \in B(\phi) $
and let $ \eta \in G(\psi), \ {\bf E} \eta = 0. $ Let also $ \xi(i), i = 1,2,\ldots $ be
independent copies of $ \xi $ and let $ \eta(j) $ be independent copies of $ \eta. $
Then

$$
\sup_n  ||n^{-1/2} \sum_{i=1}^n \xi(i)||B(\overline{\phi}) \le ||\xi||B(\phi), \eqno(3.8a),
$$
(exponential level),
$$
\sup_n  ||n^{-1/2} \sum_{j=1}^n \eta(j)||G(\overline{\psi}) \le ||\eta||G(\psi). \eqno(3.8b)
$$
(power level).\par
\vspace{2mm}
{\bf Remark 3.3.} It is important to notice that the exponential bounds for the
tail behavior for the sums of independent random variables (3.8a) does not be
obtained from the power estimations (3.8b) and conversely the power estimations (3.8b)
does not be obtained from exponential estimations (3.8a). Let us consider the
following examples, see \cite{Ostrovsky1}, p. 55-57. \par
 Let $ \{ \xi_i \}, \ i = 1,2,\ldots $ be a sequence  identical distributed centered
 r.v.  with the following tail function:

$$
U(\xi_i, u) = \exp \left(- u^r \right), \ r = \const > 0, \ u \ge 2.
$$
 Denote

 $$
 \overline{P}_r(u) =
 \sup_n {\bf P} \left(n^{-1/2} \left|\sum_{i=1}^n \xi_i \right| > u \right), \ u \ge 2.
 $$
 From the relation (3.8a) follows the estimation

 $$
 \overline{P}_r(u) \le \exp \left( - C(r) u^{\min(r,2)} \right), \ u \ge 2, \eqno(3.8c)
 $$
only under the condition $ r \ge 1. $ \par
But from the relation (3.8b) follows the inequality

$$
\overline{P}_r(u) \le \exp
\left( - C_1(r) u^{r/(r+1) } [\log u ]^{C_2(r)}  \right), \ u \ge 2. \eqno(3.8d)
$$
 Note that the inequality (3.8c) is more exact that (3.8d), but only in the case
 $ r > 1. $ \par
In the case $ r \in (0,1). $ the r.v. $ \xi_i $ does not belong to the any
$ B(\phi) $ space, $ \phi \in \Phi;$ therefore the relation (3.8d) has  advantage
in the considered variant. \par

\vspace{3mm}
{\bf Pilcrow D.} {\it Continuity and compactness of random fields.} \par

 Let $ \xi(t) = \xi(t,\omega), \ t \in T $ be a separable centered random field.
The "constructive" introduction of the $ \psi = \psi_{\xi}(p)  $ function for the
random field $ \xi(\cdot) $ may be follows:

$$
\psi_{\xi}(p) = \sup_{t \in T} |\xi(t)|_p, \eqno(3.9)
$$
if it is finite for some $ p > 1.$ \par
The natural function $ \phi_{\xi}(\lambda)  $  for the $ \xi(\cdot) $ is defined as follows:

$$
\phi_{\xi}(\lambda) = \max_{\mu = \pm 1} \sup_{t \in T}
\log {\bf E} \exp(\lambda \ \mu \ \xi(t)), \eqno(3.10)
$$
if reasonably the last function is finite for some non-trivial interval

$$
\lambda \in (-\lambda_0, \lambda_0), \ \lambda_0 = \const > 0.
$$
Evidently, $ \psi_{\xi} (\cdot) \in \Psi, \ \phi_{\xi}(\cdot) \in \Phi. $\par
 Let us denote for arbitrary function $ \psi \in \Psi $

 $$
 v_{*}(x) = v_{*,\psi}(x) \stackrel{def}{=} \inf_{y \in (0,1)} (xy + \log \psi(1/y)).
 \eqno(3.11)
$$

 M.Ledoux and M.Talagrand \cite{Talagrand1}, chapter 2
 write instead our function $ \exp \left(-\phi^*(x) \right) $ some Young's
function $ \Psi(x) $ and used as a rule a function $ \Psi(x) = \exp(-x^2/2) $
(the so-called "subgaussian case").\par

 Without loss of generality we can and will suppose

$$
\sup_{t \in T} [ \ ||\xi(t) \ ||B(\phi)] = 1,
$$
(this condition is satisfied automatically in the case of natural choosing
of the function $ \phi: \ \phi(\lambda) = \phi_0(\lambda) \ ) $
and that the metric space $ (T,d) $ relatively the so-called
{\it natural} distance (more exactly, semi-distance)

$$
d_{\phi}(t,s) \stackrel{def}{=} ||\xi(t) - \xi(s)|| B(\phi)\eqno(3.12)
$$
is complete. \par
 Recall that the semi-distance $ d = d(t,s), \ s,t \in T $ is, by definition,
non-negative symmetrical numerical function, $ d(t,t) = 0, \ t \in T, $
satisfying the triangle inequality, but the equality $ d(t,s) = 0 $
does not means (in general case) that $ s = t. $ \par
 For example, if $ \xi(t) $ is a centered Gaussian field with covariation function
 $ D(t,s) = {\bf E} \xi(t) \ \xi(s), $ then
$ \phi_0(\lambda) = 0.5 \ \lambda^2, \ \lambda \in R, $ and $ d(t,s) = $

$$
||\xi(t) - \xi(s)||B(\phi_0) = \sqrt{ \bf{Var} [ \xi(t) - \xi(s) ]} =
\sqrt{ D(t,t) - 2 D(t,s) + D(s,s) }.
$$

 Let us introduce for any subset $ V, \ V \subset T $ the so-called
{\it entropy } $ H(V, d, \epsilon) = H(V, \epsilon) $ as a logarithm
of a minimal quantity $ N(V,d, \epsilon) = N(V,\epsilon) = N $
of a balls $ S(V, t, \epsilon), \ t \in V: $
$$
S(V, t, \epsilon) \stackrel{def}{=} \{s, s \in V, \ d(s,t) \le \epsilon \},
$$
which cover the set $ V: $
$$
N = \min \{M: \exists \{t_i \}, i = 1,2,…, M, \ t_i \in V, \ V
\subset \cup_{i=1}^M S(V, t_i, \epsilon ) \},
$$
and we denote also
$$
B(t_i,\epsilon) = \log N; \ S(t_0,\epsilon) \stackrel{def}{=}
 S(T, t_0, \epsilon), \ H(d, \epsilon) \stackrel{def}{=} H(T,d,\epsilon).
$$
 A {\it capacity} $ M= M(V,d,\epsilon) $ of the set $ V $ is the maximal number of
 disjoint balls
 $$
 B(t_i,\epsilon) = \{ s, s \in V, d(t_i,s) \le \epsilon \}, \ t_i, s \in V, i = 1,2,\ldots, M,
 $$

 $$
 B(t_i,\epsilon) \cap B(t_j,\epsilon) = \emptyset, i \ne j, \ i,j = 1,2,\ldots,M,
 $$
subsets of the set $ V:$

$$
\cup_{i=1}^M B(t_i,\epsilon) \subset V.
$$
 Denote
$$
H(\epsilon) = H(T,d,\epsilon), \ L(\epsilon) = \log M(T,d,\epsilon),
$$

$$
h_-(\epsilon) = \inf_{t \in T} \mu(B(t,\epsilon)), \ h_+(\epsilon) =
\sup_{t \in T} \mu(B(t,\epsilon)).
$$
 It is known  \cite{Ostrovsky1}, p. 95-97; \cite{Vitushkin1}, p. 5-9 that
 $$
 H(2\epsilon) \le L(\epsilon) \le H(\epsilon).
 $$
  The next fact allows to estimate the values $ H(\epsilon), L(\epsilon): $

  $$
  \frac{\mu(T)}{h_+(2\epsilon)} \le M(T,\epsilon) \le \frac{\mu(T)}{h_-(\epsilon)}.
  $$
   If for instance the set $ T $  is bounded open subset of the space $ R^d $ equipped
with a distance  $ d(t,s) $   for which
$$
C_1 |t-s|^{\alpha} \le d(t,s) \le C_2 |t-s|^{\alpha}, \ \alpha = \const \in (0,1],
$$

$$
|t-s| = \sqrt{\sum_{k=1}^d (t_k-s_k)^2},
$$
then
$$
N(T,d,\epsilon) \asymp C_3(T,\alpha) \ \epsilon^{-d/\alpha}, \ \epsilon \in (0,1).
$$

 It follows from Hausdorf's theorem that
$ \forall \epsilon > 0 \ \Rightarrow H(V,d,\epsilon)< \infty $ iff the
metric space $ (V, d) $ is precompact set, i.e. is the bounded set with
compact closure.\par
\vspace{3mm}
 Let $ \xi(t), \ t \in T $  be a separable numerical random field
 such that  for some function $ \psi \in \Psi \ \sup_{t \in T } ||\xi(t)||G(\psi) =1. $
  We introduce so-called {\it natural distance}

  $$
  d_{\psi}(t,s) \stackrel{def}{=} ||\xi(t) - \xi(s)||G(\psi). \eqno(3.13)
  $$
 Denote

 $$
 v(y) = \log \psi(1/y), \  \ v_*(x) = \inf_{y: \psi(1/y) < \infty} (xy + v(y)).
 $$

{\bf Lemma 3.1.a}  {\it If the following integral converges:}

$$
K(\sigma_{\phi}):= \int_0^{\sigma} \chi_{\phi}(H(T,d_{\phi},x)) \ dx < \infty,
$$
{\it where }
$$
\sigma_{\phi} = \sup_{t \in T} ||\xi(t)||B(\phi) < \infty,
$$

{\it then (see \cite{Kozachenko1},
\cite{Ostrovsky1}, chapter 3, section 3.4)}

$$
{\bf P} \left(\xi(\cdot) \in C(T,d_{\phi} \right)) = 1
$$
{\it and for the values}

$$
u \ge 2 K \left(\sigma_{\phi} \right)
$$
{\it we have}
$$
{\bf P} \left( \sup_{t \in T} |\xi(t)| > \sigma_{\phi} u \right) \le
2 \exp \left(- \phi^* \left( u - \sqrt{2 K(\sigma_{\phi}) u} \right)  \right).\eqno(3.14a)
$$

\vspace{3mm}

{\bf Lemma 3.1.b}  {\it If the following integral converges:}

   $$
   \int_0^1 \exp \left(v_*(\log (2 H(T, d_{\psi}, x))) \right) \ dx  < \infty,
   $$
{\it then (see \cite{Ostrovsky1}, chapter 4, section 4.1)}

$$
{\bf P} \left(\xi(\cdot) \in C(T,d_{\psi} \right)) = 1.
$$
{\it Moreover,}

 $$
 ||\sup_{t,s: d_{\psi}(t,s) \le \delta } |\xi(t) - \xi(s)| \ ||G(\psi) \le
 Z(\delta),
$$
where
$$
Z(\delta) = Z(\psi,\delta) \stackrel{def}{=}
 9 \int_0^{\delta} \exp \left(v_*(\log (2 H(T, d_{\psi}, x))) \right) \ dx.
 \eqno(3.14b)
 $$
 As a slight consequence: let us denote
$$
\sigma_{\psi} = \sup_{t \in T} ||\xi(t)||G(\psi) < \infty.
$$
  We have for arbitrary fixed non-random value $ t_0 \in T $ using triangle
inequality  and taking in (3.14b) $ \delta = \sigma_{\psi}:  $

$$
\sup_{t \in T} |\xi(t)| \le |\xi(t_0)| + \sup_{t \in T}|\xi(t) - \xi(t_0)|;
$$

$$
||\sup_{t \in T} |\xi(t)| \ ||G(\psi) \le ||\xi(t_0)||G(\psi) +
Z \left(\sigma_{\psi} \right) \le
$$

 $$
 \sigma_{\psi} + Z \left(\sigma_{\psi} \right) \stackrel{def}{=}
 \overline{Z} = \overline{Z}(\psi) = \overline{Z}(\psi, \epsilon);
 $$

$$
{\bf P} \left(\sup_{t \in T} |\xi(t)| > u \right)  \le
\exp \left[- w^*_{\psi} \left(u/\overline{Z}(\psi) \right) \right], \
u \ge 2 \overline{Z}(\psi). \eqno(3.14c)
$$

\vspace{3mm}

   For instance, let $ \psi(p) = p^{\beta}, \beta = \const > 0.  $ The condition (3.14b) may be written as

 $$
 \int_0^1 H^{\beta}(T,d_{\psi}, x) dx < \infty.
 $$
 Moreover, if

 $$
 \psi(p) = p^{\beta_1} \ \log^{\beta_2}(1+p),  \ p \ge 1, \ \beta_1 = \const > 0,
 $$
then the condition (3.14) has a view

 $$
 \int_0^1 H^{\beta_1}(T,d_{\psi}, x) \
 \log^{\beta_2} \left(1+ H(T,d_{\psi},x) \right) \ dx < \infty.
 $$
 Let now $ \supp \psi = [1,b], \ b = \const > 1. $ The condition (3.14) is equivalent to the famous Pizier's condition \cite{Pizier1}

 $$
 \int_0^1 N^{1/b}(T, d_{\psi}, x) \ dx < \infty.
 $$

\vspace{3mm}
 Further, let $ \xi_{\alpha}(t), \ t \in T, \ \alpha \in  A $  be arbitrary
 {\it family } of random fields,  where $ A $  is arbitrary set, such that
such that  for some function $ \psi \in \Psi $

$$
\sup_{\alpha \in A} \sup_{t \in T } ||\xi_{\alpha}(t)||G(\psi) =1.
$$
  We introduce again the so-called {\it natural distance} induced by the
  family $ \{ \xi_{\alpha}(\cdot)  \} $

  $$
  d_{\psi}(t,s) \stackrel{def}{=} \sup_{\alpha \in A} ||\xi_{\alpha}(t) -
  \xi_{\alpha}(s)||G(\psi). \eqno(3.15)
  $$
 Denote

 $$
 v(y) = \log \psi(1/y), \  \ v_*(x) = \inf_{y: \psi(1/y) < \infty} (xy + v(y)).
 $$
{\bf Lemma 3.2.}
{\it Assume that for some $ t_0 \in T $ the one-dimensional family of random  variables
 $ \{\xi_{\alpha} (t_0)   \} $ is stochastically bounded:}

 $$
 \lim_{u \to \infty} \sup_{\alpha \in A} {\bf P}(|\xi_{\alpha}(t_0)| > u) = 0. \eqno(3.16)
 $$

{\it   If the following integral converges:}

   $$
   \int_0^1 \exp \left(v_*(\log (2 H(T, d_{\psi}, x))) \right) \ dx  < \infty,
   \eqno(3.17)
   $$
{\it then  the family of distributions $ \mu_{\alpha}(\cdot) $ generated by the random fields $ \xi_{\alpha}(\cdot) $ in the space} $ C(T, d): $

$$
\mu_{\alpha}(B) = {\bf P} \left(\xi_{\alpha}(\cdot) \in B \right)
$$
{\it is weakly compact.}\par

 The detail explanation of the theory of weak convergence for probabilistic measures in
the metric spaces see in the monographs \cite{Prokhorov1},
 \cite{Billingsley1}, \cite{Pollard1}.  A main
conclusion of this theory for the continuous random processes may be formulated as
follows. If the finite-dimensional distributions of the sequence of a random processes
$ \eta_n(t), \ n=1,2,\ldots $ converge as $ n \to \infty $ to the finite-dimensional
distributions of a random process $ \eta(t) $ and the distributions of the random processes
$$
\mu_n(B) = {\bf P}(\eta_n(\cdot) \in B)
$$
are weakly compact, then the distributions of arbitrary continuous functional
$ F: C(T) \to R $ at the points $ \eta_n(\cdot) $ converge to the distribution
$  F(\eta): \ \forall \lambda \in R $
$$
\lim_{n \to \infty} {\bf E}\exp(i \lambda F(\eta_n))= {\bf E}\exp(i \lambda F(\eta)).
$$
In particular:

$$
\lim_{n \to \infty} {\bf P}( \sup_{t \in T}|\eta_n(t)| > u )=
{\bf P}( \sup_{t \in T}|\eta(t)| > u ) \eqno(3.18)
$$
for all the points $ u, \ u \in R $ in which the function
$$
u \to {\bf P}( \sup_{t \in T}|\eta(t)| > u )
$$
is continuous. \par
\vspace{2mm}
{\bf Pilcrow E.} {\it Central Limit Theorem (CLT) in Banach space.} \par
\vspace{2mm}
Let $ (\Omega, {\cal A}, {\bf P}) $ be a probabilistic space and let $ {\cal B} $
with a norm $ |||\cdot||| $  be a Banach space. We will say as ordinary that the (centered) random variable $ \zeta, \ \zeta: \Omega \to {\cal B} $ satisfies the CLT
in this space, if the sequence
$$
\overline{\zeta}_n = n^{-1/2}  \sum_{i=1}^n \zeta_i \eqno(3.19)
$$
converges weakly in distribution as $ n \to \infty $ to non-trivial Gaussian centered
random variable $ \overline{\zeta}_{\infty}: $

$$
\lim_{n \to \infty} \Law(\overline{\zeta}_n) = \Law(\overline{\zeta}_{\infty}). \eqno(3.20)
$$
 Here the variables $ \{ \zeta_i \} $ are independent copies $ \zeta. $ \par
 It is evident that the variable $ \overline{\zeta}_{\infty} $ has at the same
 covariation operator as the variable $ \zeta. $
 If (3.20) there holds, then for all positive values $ u $

$$
\lim_{n \to \infty} {\bf P} (|||\overline{\zeta}_n||| > u) =
{\bf P} (|||\overline{\zeta}_{\infty}||| > u). \eqno(3.21)
$$
As well as in the one-dimensional case,  the CLT in the Banach space will be used
in the Monte-Carlo  method for building of confidence region for estimated function
in the Banach space norm.\par

\vspace{2mm}
{\bf Pilcrow F.} {\it Multiplicative inequalities in Grand Lebesgue spaces.} \par
\vspace{2mm}

 Let $ \xi \in G(\psi_1), \eta \in G(\psi_2), \ \tau = \xi \cdot \eta. $ We will study
 in this subsection the inequalities of a view

 $$
||\tau||G(\psi_3) = ||\xi \eta||G(\psi_3) \le C ||\xi||G(\psi_1) \ ||\eta||G(\psi_2).
\eqno(3.22)
 $$
("power" level) or analogously

 $$
||\tau||B(\phi_3) = ||\xi \eta||B(\phi_3) \le C ||\xi||B(\phi_1) \ ||\eta||G(\phi_2).
\eqno(3.23)
 $$
("exponential" level). \par
It is convenient to continue arbitrary function  $  \psi = \psi(p) $ as follows:

$$
 \forall p \notin \supp( \psi) \Rightarrow \psi(p) = +\infty.\eqno(3.24)
$$

\vspace{2mm}
{\bf A. Independent case, power level.} \par
\vspace{2mm}
 Note first of all that if the r.v. $ \xi, \eta  $ are independent, then

 $$
 |\xi \eta|_p = |\xi|_p \ |\eta|_p.
 $$
 Therefore, if $ \xi \in G(\psi_1), \ \eta \in G(\psi_2), $ then \par
{\bf Proposition 3.A.}
 $$
||\xi \eta||G(\psi_1 \cdot \psi_2) \le  ||\xi||G(\psi_1)  \cdot ||\eta||G(\psi_2).
\eqno(3.25)
 $$

\vspace{2mm}
{\bf B. Dependent case, power level.} \par
\vspace{2mm}

 We do not suppose in this pilcrow the r.v. $ \xi, \eta $ to be independent. Let again
$ \xi \in G(\psi_1), \ \eta \in G(\psi_2). $  We introduce the following operation
for two functions $ \psi_1(\cdot), \ \psi_2(\cdot) $  from the set  $ \Psi: $

$$
\psi_1 \circ \psi_2 (r) \stackrel{def}{=} \inf_{p > 1}
\{\psi_1(pr) \cdot \psi_2(rp/(p-1)) \}. \eqno(3.26)
$$
 If $ \supp \psi_1 = [1,b_1), \supp \psi_2 = [1,b_2), $ then

 $$
 \supp \left(\psi_1 \circ \psi_2 (\cdot) \right) = [1, b_3),
 $$
 where
 $$
b_3 = b_1 b_2/(b_1 + b_2). \eqno(3.27).
 $$
 Note that if  $ (b_1-1)(b_2-1) > 1, $ then $ b_3 > 2 $ and if
$ (b_1-2)(b_2-2) > 4, $ then $ b_3 > 2. $ \par
 It is known \cite{Liflyand1} that in the general, i.e. dependent case \par
{\bf Proposition 3.B.}
$$
||\xi \ \eta||G(\psi \circ \psi_2) \le ||\xi||G(\psi_1)  \cdot ||\eta||G(\psi_2).
\eqno(3.28)
$$

\vspace{2mm}
{\bf C. Independent case, exponential level.} \par
\vspace{2mm}

 Let $ \xi \in B(\phi_1) $ and $ \eta \in B(\phi_2)$ be independent r.v. and
 $  \ \phi_{1,2}(\cdot) $ be two functions from the set $  \Phi $ such that

$$
\overline{\lim}_{p \to \infty}
\frac{\log [p^2/(\phi_1^{-1}(p) \phi_2^{-1}(p))]}{\log p} \le 1. \eqno(3.29)
$$
 We define a new operation (commutative and associative)
 $ \phi_3(p) = \phi_1 \odot \phi_2(p) $ for the functions
 $ \phi_{1}(\cdot) $ and $ \phi_{2}(\cdot) $ as follows:

$$
\phi_3(p) = \left[\phi_1^{-1}(p) \ \phi_2^{-1}(p)/p  \right]^{-1}, \ p \ge 1. \eqno(3.30)
$$

 Recall that at $ |\lambda| \le 1 \ \Rightarrow \phi_3(\lambda) = C \cdot \lambda^2,  $
where the constant $ C $ must be choose such that  the function $ \lambda \to
\phi_3(\lambda) $  is continuous. \par

Note that by virtue of independence

$$
|\xi \ \eta|_p \le C_2^2 \ \frac{p}{\phi_1^{-1}(p)} \ \frac{p}{\phi_2^{-1}(p)} \le
C_3 \frac{p}{\phi_3^{-1}(p)},
$$
therefore under condition  (3.29) \par
{\bf Proposition 3.C.}
$$
||\xi \ \eta||B(\psi_3) \le C_4 \ ||\xi||B(\phi_1) \cdot ||\eta||B(\phi_2). \eqno(3.31)
$$
 As a consequence: if the r.v. $ \xi_1, \xi_2, \ldots, \xi_k $ are mutually  independent and $ \xi_j \in B(\phi_j),  $ then

 $$
 ||\xi_1 \ \xi_2 \ \ldots \xi_k||B(\phi^{(k)}) \le C^k_4 ||\xi_1||B(\phi_1) \
 ||\xi_2||B(\phi_2) \ldots ||\xi_k||B(\phi_k), \eqno(3.32)
 $$
where

$$
\phi^{(k)} = (((\phi_1 \odot \phi_1)  \odot \phi_2) \ldots \odot \phi_k) . \eqno(3.33)
$$
 If for example the r.v. $ \xi, \ \eta $  are subgaussian and independent:

 $$
 ||\xi||\sub = \sigma_1, \ ||\eta||\sub = \sigma_2, \ \sigma_i < \infty,
 $$
then

$$
|\xi|_p \le C_4 \sigma_1 \ \sqrt{p}, \ |\eta|_p \le C_4 \sigma_2 \ \sqrt{p},
$$
and the r.v. $ \tau = \xi \cdot \eta $ satisfies the following moment condition:

$$
|\tau|_p \le C_5 \ \sigma_1 \ \sigma_2 \ p, \ p \in [1,\infty),
$$
or equally the r.v. $ \tau $ satisfies the Kramer's condition. \par

\vspace{2mm}
{\bf D. Dependent case, exponential level.} \par
\vspace{2mm}

Let again $ \xi \in B(\phi_1) $ and $ \eta \in B(\phi_2)$ be arbitrary r.v. and
 $  \ \phi_{1,2}(\cdot) $ be two functions from the set $ \Phi. $
 We define a new operation (commutative and associative)
 $ \phi_4(p) = \phi_1 \otimes \phi_2(p) $ for the functions
 $ \phi_{1}(\cdot) $ and $ \phi_{2}(\cdot) $ as follows. Denote

 $$
 \psi_j(p) = \frac{p}{\phi^{(-1)}_j(p)}, \ j = 1,2;
 $$

$$
\psi_4(p) = [\psi_1 \circ \psi_2](p),
$$
then we define

$$
\phi_4(p) = [\phi_1 \otimes \phi_2](p) \stackrel{def}{=}
\left[ \frac{p}{\psi_4(p)} \right]^{(-1)}, \eqno(3.34)
$$
if obviously for some
$$
b > 1 \Rightarrow \phi_4(b) < \infty.  \eqno(3.35)
$$
 Note that for the CLT in the Banach space we need to assume $ b \ge 2. $\par
The function $ \phi_4(p) = [\phi_1 \otimes \phi_2](p) $ has the following sense.
If the condition (3.35) is satisfied, then \par
{\bf Proposition 3.D.}
$$
||\xi \ \eta ||B(\phi_4) \le C \ ||\xi||B(\phi_1) \cdot ||\eta||B(\phi_2). \eqno(3.36)
$$

 As a consequence: if the r.v. $ \xi_1, \xi_2, \ldots, \xi_k $ are mutually independent and $ \xi_j \in B(\phi_j),  $ then

 $$
 ||\xi_1 \ \xi_2 \ \ldots \xi_k||B(\phi_{(k)}) \le C^k_4 ||\xi_1||B(\phi_1) \
 ||\xi_2||B(\phi_2) \ldots ||\xi_k||B(\phi_k), \eqno(3.37)
 $$
where

$$
\phi_{(k)} = (((\phi_1 \otimes \phi_1)  \otimes \phi_2) \ldots \otimes \phi_k). \eqno(3.38)
$$

 Let us return for instance to the subgaussian case.
 If  the r.v. $ \xi, \ \eta $  are subgaussian (and arbitrary dependent):

 $$
 ||\xi||\sub = \sigma_1, \ ||\eta||\sub = \sigma_2, \ \sigma_i < \infty,
 $$
then the centered r.v. $ \tau = \xi \cdot \eta - {\bf E} (\xi \eta) $
satisfies the following moment condition:

$$
|\tau|_p \le C_6 \ \sigma_1 \ \sigma_2 \ p, \ p \in [1,\infty), \ C_6 > C_5,
$$
or again the r.v. $ \tau $ satisfies the Kramer's condition. \par
 Let us consider more essential but more exotic example. Let $ \{\xi \} $  be
 symmetrically distributed r.v. with the following tail behavior:

$$
{\bf P}(\xi > x) \le \exp \left( - C_1 \ e^x  \right), \  x > 1, \eqno(3.39)
$$

and put

$$
\upsilon = \upsilon_k = \prod_{i=1}^k \xi_i - {\bf E}\prod_{i=1}^k \xi_i,
$$
where $ \xi_i  $ are independent copies of $ \xi. $ It follows from proposition 3.D
after some calculations

$$
U(\upsilon_k,x) \le \exp \left( - C_2(k,C_1) \ e^{x^{1/k}}  \right), \  x > 1, \eqno(3.40)
$$
 In order to illustrate the inequality (3.40), let us consider the following example.
Let  $ \eta_i = \zeta, \ i = 1,2,\ldots, \ \zeta > 0, $ and

$$
{\bf P}(\zeta > x) = \exp \left( - \ e^x  \right), \  x > 1,
$$
and  let $ \nu = \prod_{i=1}^k \eta_i = |\zeta|^k; $  then

$$
{\bf P}(\nu > x) = {\bf P} \left(|\zeta| > x^{1/k} \right) =
\exp \left( - e^{x^{1/k}}  \right), \  x > 1.
$$

\vspace{3mm}
{\bf Remark 3.D.} Let us consider the $ B \left(\phi_{e,\kappa} \right),
\kappa \ge 1 $ space as a $ B(\phi) $ space with the correspondent function

$$
\phi_{e,\kappa}(\lambda) \asymp |\lambda| \cdot \log^{\kappa}(2 + |\lambda|),
|\lambda| \ge 1,  \eqno(3.41)
$$
or equally

$$
\log [\phi^*_{e,\kappa}(x)] \asymp \left( C(\kappa) x^{1/\kappa} \right),  \ x > 1,
$$

$$
\psi_{\phi_{e,\kappa}}(p) \asymp \log^{\kappa}(p+1), \ p \ge 1.
$$
 Let also $ \xi = \xi(t), \ t \in T $ be separable random field such that

$$
\sup_{t \in T} ||\xi(t)||B \left(\phi_{e,\kappa} \right) = 1.
$$
 Introduce the following (finite) distance

 $$
 \rho_{e,\kappa}(t,s) =  ||\xi(t) - \xi(s)||B \left(\phi_{e,\kappa} \right).
 \eqno(3.42)
 $$

If the following integral converges:

$$
I_{e,\kappa}= \int_0^1 \log^{2 \kappa} \left( 1 + H(T,\rho_{e,\kappa},x) \right) \ dx < \infty, \eqno(3.43)
$$
then
$$
{\bf P}(\xi(\cdot) \in C(T,\rho_{e,\kappa})) = 1 \eqno(3.44)
$$
and

$$
{\bf P}(\sup_{t \in T} |\xi(t)| > u) \le
\exp \left( - C_3 \exp \left[C_4(I_{e,\kappa}) x^{1/\kappa} \right] \right), \ x \ge 1.\eqno(3.45)
$$

\vspace{2mm}
{\bf Pilcrow G.} {\it Dual spaces for Grand Lebesgue spaces.} \par
\vspace{2mm}
 We describe briefly in this pilcrow the {\it dual (conjugate)} spaces to  the
 Grand Lebesgue Spaces. \par
 We do not need to suppose (only in this Pilcrow G!) the finiteness condition
 $ \mu(T) = 1; $ it is sufficient to entrust on the
measure $ \mu $ instead the condition  of sigma-finiteness; and suppose also
the triplet $ (T,\Sigma, \mu) $ is resonant in the terminology of the classical book \cite{Bennet1}. This imply by definition that either the measure $ \mu $ is diffuse:

$$
\forall A \in \Sigma, \ 0 < \mu(A) < \infty \Rightarrow  \exists B \subset A,  \
\mu(B) = \mu(A)/2
$$
or the measure $ \mu $ is purely discrete and each atoms have at the same (positive)
weight.
\par
  Further, let $ (a,b), 1 \le a < b \le \infty  $ be a maximal open subset of the
 support of some  function $ \psi(\cdot) \in \Psi. $ We will consider here only non-trivial case when

 $$
 \max(\psi(a+0), \psi(b-0)) = \infty. \eqno(3.45)
 $$

 Note that the associate space for GLS spaces are describes in \cite{Ostrovsky20};
a less general case see in \cite{Fiorenza2}, \cite{Fiorenza3}.\par
 Namely, let us introduce the space $ DG(\psi) = DG\psi(a,b) $ consisting on all the measurable functions $  \{g \}, \ g: T \to R  $ with finite norm

 $$
 ||g||DG(\psi)= \inf_{ \{p(k), p(k) \in (a,b) \}} \inf_{\{g_k \}}
 \{ \sum_k \psi(p(k)) |g_k|_{p(k)/(p(k)-1)}  \},\eqno(3.46)
 $$
where interior $ \inf $  is calculated over all the sequences $ \{ g_j \}, $  finite
or not, of a measurable functions such that

$$
g(x) = \sum_j g_j(x),
$$
and exterior  $ \inf $  is calculated over all the sequences $ \{ p(k) \}, $ belonging
to the {\it open} interval $ (a,b). $\par
 An action of a linear continuous functional $ l_g, \ g \in DG(\psi) $  on the
 arbitrary function $ f \in G(\psi) $ may be described as ordinary by the formula

 $$
 l_g(f) = \int_T f(x) \ g(x) \ \mu(dx)
 $$
 with $ ||l_g|| = ||g||DG(\psi). $ \par
  Recall that the set of all such a functionals $ \{ l_g \} $ equipped with the norm
$ ||l_g|| = ||g||DG(\psi) $ is said to be associate space to the space $ G(\psi) $
and is denoted as usually $ G'(\psi) = [G(\psi)]'.$ \par

  Define the space $ G^o(\psi) $ as a (closed) subspace of a space $ G(\psi) $
 consisting on all the functions $ \{g\} $ from the space $ G(\psi) $  satisfies
 the condition

 $$
 \lim_{\psi(p)\to \infty} \frac{|g|_p}{\psi(p)} = 0, \eqno(3.47)
 $$
and introduce correspondent quotient space

$$
G_o(\psi) = G(\psi)/G^o(\psi). \eqno(3.48)
$$
 It is known  \cite{Ostrovsky20} that the space $ DG(\psi) $ is associate to the space
 $ G(\psi) $ and is dual to the space $ G^o(\psi). $ \par
It is easy to verify analogously to the case of Orlicz space (see \cite{Rao2}, p. 119-121, \cite{Rao3}, \cite{Morse1}, \cite{Luxemburg1} that each  function $ f \in G(\psi) $
may be uniquely represented as a sum

$$
f = f^o + f_o,  \eqno(3.49)
$$
(direct sum), where $f^o \in G^o(\psi), \ f_o \in G_o(\psi). $ \par
 More exactly, the function $ f_o $ represented the class of equivalence under
 relation

 $$
 f_1 \sim f_2 \ \Leftrightarrow f_1 - f_2 \in G^o(\psi).
 $$
Therefore, arbitrary continuous linear functional $ L = L_{g,\nu}(f) $ on the space $ G(\psi) $ may be uniquely represented as follows:

$$
L_{g,\nu}(f) = \int_T f^o(x) \ g(x) \ \mu(dx) + \int_T f_o(x) \ \nu(dx), \eqno(3.50)
$$
where $ g \in DG(\psi), \ \nu \in ba(T,\Sigma,\mu); \ ba(T,\Sigma,\mu) $ denotes
the set of all {\it finite additive} set function with finite total variation:

$$
|\nu| =|\nu|(T) = \sup_{A \in \Sigma} [\nu(A) - \nu(T \setminus A)] < \infty. \eqno(3.51)
$$
 Note that the (generalized) measure $ \nu $ is singular relative the source measure
 $  \mu; $  therefore

 $$
 ||L_{g,\nu}||G^*(\psi) = ||g||DG(\psi) + |\nu|. \eqno(3.52)
 $$

\vspace{2mm}
\section{ Monte-Carlo method for the parametric integrals calculation}
\vspace{2mm}
 We consider in this section the problem of Monte-Carlo approximation and construction
 of a confidence region  in the uniform norm for the parametric integral of a view

 $$
 I(t) = \int_X g(t, x) \ \nu(dx). \eqno(4.1)
 $$
 Here $ (X, \Sigma, \nu) $ is also a probabilistic space with normed:
$ \nu(X) = 1 $ non-trivial measure $ \nu. $ \par
 A so-called "Depending Trial Method" estimation for the integral (4.1) was introduced by
Frolov A.S.and Tchentzov N.N., see \cite{Frolov1}:

$$
I_n(t) = n^{-1} \sum_{i=1}^n g(t,\eta_i), \eqno(4.2)
$$
where $ \{\eta_i\} $ is the sequence of $ \nu $ distributed:

$$
{\bf P} (\eta_i \in A) = \nu(A)
$$
independent random variables.\par
 We intend in this section to improve the result of the article \cite{Frolov1} and
its consequence,  see \cite{Ostrovsky1}, chapter 5, section 5.11. \par
The modern methods of (pseudo)random variable generations are described,
in particular, in \cite{Goldreich1}; see also \cite{Devroye1}.\par
We assume $ \forall t \in T \ g(t,\cdot) \in L_2(X,\nu): $

$$
\int_X g^2(t,x) \ \nu(dx) < \infty;
$$
then the integral $ I(t) $ there exists for all the values $ t; t \in T $ and we can for
any {\it fixed} point $ t_0 \in T $ use for an error evaluating the classical Central
Limit Theorem:

$$
\lim_{n \to \infty} {\bf P}( \sqrt{n} |I_n(t_0) - I(t_0)| \le u )=
\Phi(u/\sigma_0) - \Phi(-u/\sigma_0), \eqno(4.3)
$$
where as usually

$$
\Phi(u) = (2 \pi)^{-1/2} \int_{-\infty}^u \exp(-z^2/2) \ dz
$$
and

$$
\sigma^2(t) = {\bf Var} (f(t,\eta)) = \int_X g^2(t,x) \nu(dx) - I^2(t),
$$
$ \sigma_0 = \sigma(t_0). $
 Let us consider now the problem of building confidence region for $ I(t) $ in the
 uniform norm, i.e. we investigate the probability

 $$
 {\bf P^{(n)}}(u) \stackrel{def}{=}
 {\bf P}(\sup_{t \in T}\sqrt{n}| I_n(t) - I(t)| > u), \ u = \const > 0.\eqno(4.4)
 $$
 On the other words, we use for construction of confidence region in the uniform norm
 in the parametrical case the Central Limit Theorem (CLT)
 in the space of continuous functions $ C(T) $  alike in the classical case of
 ordinary Monte-Carlo method is used customary CLT.\par

\vspace{2mm}

Let us introduce as in the first section the following function, presumed to be finite
$ \nu \ - $ almost everywhere:

$$
Q(x) = \vraisup_{t \in T} |g(t,x)|  \eqno(4.5)
$$
and we introduce also the so-called a new {\it natural } distance, more exactly, semi-distance, $ \beta = \beta(t,s) $ on the space $ T: $

$$
\beta(t,s) \stackrel{def}{=} \vraisup_{x \in T} \frac{|g(t,x)-g(s,x)|}{Q(x)}, \eqno(4.6)
$$
so that

$$
|g(t,x) - g(s,x)| \le Q(x) \ \beta(t,s). \eqno(4.7)
$$
 Let us consider the centered random process

 $$
 g^0(t,x) = g(t,x) - I(t); \ {\bf E} g^0(t,\eta)=0,
 $$
and define the $ \psi - $ functions as follows:

$$
\psi(p) = \psi_g(p) = \sup_{t \in T} |g^0(t,\eta)|_p = \sup_{t \in T}
{\bf E^{1/p}} |g^0(t,\eta)|^p =
$$

$$
\left[\int_X |g(t,x) - I(t)|^p \ \nu(dx)  \right]^{1/p}. \eqno(4.8)
$$

 We suppose the function $ \psi(\cdot) $ is finite on some interval $ 2 \le p \le b,  $
where $ b = \const \le \infty. $
 It is evident that in the case $ b  < \infty \ \overline{\psi}(p) \asymp \psi(p). $ \par
We introduce on the basis of the function $ \overline{\psi}(p) $ the new distance
on the set $ T: $

$$
\gamma_{\psi}(t,s) = || g^0(t,\eta) - g^0(s,\eta)||G(\psi_g).\eqno(4.9)
$$

 Further, put

 $$
 \phi(\lambda) = \phi_g(\lambda) = \max_{\varepsilon = \pm 1} \sup_{t \in T} \log {\bf E} \left( \exp( \varepsilon \lambda g^0(t,\eta))\right) =
 $$

$$
\max_{\varepsilon = \pm 1} \sup_{t \in T} \log \int_T
\exp(\varepsilon \ \lambda \ g^0(t,x)) \ \nu(dx). \eqno(4.10)
$$

 If  the function $ \phi(\cdot) $ is finite on some non-trivial interval
 $ \lambda \in (-\lambda_0, \lambda_0),  $ where $ \lambda_0 = \const > 0, $
 we introduce on the basis of the function
 $$
 \overline{\phi_g}(\lambda) = \sup_{n=1,2,\ldots} n \ \phi_g(\lambda/\sqrt{n})
 $$
 the new distance on the set $ T: $

$$
\gamma_{\phi}(t,s) = || g^0(t,\eta) - g^0(s,\eta)||B(\phi_g).\eqno(4.11)
$$
 Obviously,

$$
  \gamma_{\psi}(t,s) \le C_1 \beta(t,s), \ \gamma_{\phi}(t,s) \le C_2 \ \beta(t,s).\eqno(4.12)
$$

 We denote for arbitrary {\it separable} numerical bounded with probability one
 random field  $ \zeta(t) = \zeta(t,\omega); \ t \in T, $ where
 $ ( \Omega, \cal{B}, {\bf P} ) $ is probabilistic space,

$$
{\bf P}_{\zeta}(u) = {\bf P} \left(\sup_{t \in T} |\zeta(t)| > u \right), \eqno(4.13)
$$
and correspondingly

$$
{\bf P^{(n)}}(u):= {\bf P} \left(\sqrt{n} \sup_{t \in T} |I_n(t) - I(t)| > u \right).
$$

We introduce also the centered separable Gaussian random process, more exactly,
random field $ X(t) = X_g(t), \ t \in T $ with a following covariation function

$$
Z(t,s) = Z_g(t,s) = \cov(X(t),X(s)) = {\bf E} X(t) X(s) =
$$
$$
 \int_X g(t,x) g(s,x) \ \nu(dx) - I(t)I(s). \eqno(4.14)
$$

\vspace{3mm}

{\bf Theorem 4.1. Exponential level.} \par

\vspace{2mm}

Suppose the following condition is satisfied:

$$
\int_0^1 v_{*,\overline{\psi_{\phi}}}(H(T, \gamma_{\phi}, \epsilon)) \ d \epsilon < \infty. \eqno(4.15)
$$
Then the centered Gaussian field $ X(t) $ is continuous a.e. relative the distance
$ \gamma_{\phi} $ and

$$
\lim_{n \to \infty} {\bf P^{(n)}}(u) = {\bf P}_X(u). \eqno(4.16)
$$

\vspace{3mm}

{\bf Theorem 4.2. Power level.} \par

\vspace{2mm}

Suppose the following integral converges:

$$
\int_0^1 v_{*,\overline{\psi}}(H(T, \gamma_{\psi}, \epsilon)) \ d \epsilon < \infty. \eqno(4.17)
$$
Then the centered Gaussian field $ X(t) $ is continuous a.e. relative the distance
$ \gamma_{\psi} $ and there holds

$$
\lim_{n \to \infty} {\bf P^{(n)}}(u) = {\bf P}_X(u). \eqno(4.19)
$$

\vspace{3mm}
{\bf Remark 4.1.} The exact asymptotic for the probability $ {\bf P}_X(u) $ as
$ u \to \infty $ is obtained in \cite{Piterbarg1}, p. 19, 88, 106, 114, 180:

$$
{\bf P}_X(u) \sim  C(X,T) \ u^{\kappa - 1} \exp \left(-u^2/(2\sigma_+^2) \right),
$$
where  $ C(X,T) = \const \in (0,\infty), $
$$
\sigma_+^2 = \max_{t \in T} Z_g(t,t) = \max_{t \in T}
\left[ \int_T g^2(t,x) \nu(dx) - I^2(t) \right],
$$
 the value $ \kappa = \const $ dependent on the geometrical characteristic of the set

 $$
 T_0 = \{t,s: Z_g(t,s) \in [\sigma_+^2/2, \sigma_+^2] \}.
 $$
 A non-asymptotical estimation of $ {\bf P_X(u)} $ for $ u \ge 2 \sigma_+ $ of a view

$$
{\bf P}_X(u) \le  C^+(X,T) \ u^{\kappa - 1} \exp \left(-u^2/(2\sigma_+^2) \right),
$$
 is obtained, e.g., in \cite{Ostrovsky1}, chapter 4, section 4.9. \par
\vspace{2mm}
{\bf Remark 4.2.} It is important  for the practical using of offered method to
calculate the main  parameters $ \sigma_+^2 $ and $ \kappa.  $ It may be implemented
by the following method, again by means of Monte-Carlo method. \par
 Let $ \tilde{T} = \{ t_m  \} \subset T $ be some finite net on the whole set $  T. $
 The consistent estimation $ \hat{\sigma}_+^2  $ of a value $ \sigma_+^2  $ has
 a view:

 $$
\hat{\sigma}_+^2 \approx \max_{t_m \in \tilde{T} } Z_g(t_m,t_m) \approx \max_{t_m \in \tilde{T}} \left[ n^{-1} \sum_{i=1}^n g^2(t_m,\xi_i)- I_n^2(t_m) \right].
 $$
 Analogously may be computed the distances $ d_{\phi}, \ d_{\psi}(t,s) $ and $ D(t,s); $
 for example, $ d(t_l,t_m) \approx \hat{d}(t_l.t_s), $ where

$$
\hat{d}_{\phi}(t_l,t_m) = ||g(t_l,\cdot) - g(t_m, \cdot)||B(\phi).
$$
 The consistent estimation of a value $ ||\eta||B(\phi) $ based on the independent
sample $ \{ \eta_i \}, \ i = 1,2,\ldots,n $ is described, e.g., in the monograph
\cite{Ostrovsky1}, p. 291-294. \par

  This facts may be used by building of confidence region in the uniform norm for the
 integral $ I(t). $ Namely, let $ \delta  $ be a "small"  number, for example,
 $ \delta = 0.05 $  or $ \delta = 0.01 $ etc.\par
  The value $ 1 - \delta $ may be interpreted as a reliability of confidence region.\par
  We define the value $ u(\delta) $ as a maximal solution of an equation

 $$
 {\bf P}_X(u(\delta)) = \delta,
 $$
 or asymptotically equivalently, the maximal positive solution of an equation

 $$
C(X,T) \ u(\delta)^{\kappa - 1} \exp \left(-u(\delta)^2/(2\sigma_+^2) \right) = \delta.
 $$
 The asymptotical as $ \delta \to 0+ $  confidence interval in the uniform norm
 with reliability (approximately) $ 1 - \delta $ for $ I(\cdot) $ has a view

 $$
 \sup_{t \in T} |I(t) - I_n(t)| \le \frac{u(\delta)}{\sqrt{n}}. \eqno(4.20)
 $$

\vspace{3mm}

{\bf Proof} of both theorems 4.1.and 4.2. \par
{\bf 1.} From the classical Central Limit Theorem (CLT) for the independent identically
distributed centered random vectors follows that the finite-dimensional distributions of the random fields

$$
X_n(t) = \sqrt{n} \left( I_n(t) - I(t)  \right) \eqno(4.21)
$$
converge in distribution as $ n \to \infty $ to the finite-dimensional distributions of
Gaussian field $ X(t). $  It remains to prove the {\it weak compactness} of the set
of (probabilistic) measures in the Banach space of continuous functions
$ C(T, \gamma_{\phi}) $ or correspondingly $ C(T, \gamma_{\psi}) $ induced by the
random fields $ X_n(\cdot). $ \par
{\bf 2.} We use further the result of Lemma 3.2. Namely, we put $ A = 1,2,\ldots $ and
consider the differences

$$
[I_n(t)-I(t)] - [I_n(s) - I(s)] = n^{-1/2}
\sum_{i=1}^n [(g(t,\xi_i) - I(t)) - (g(s,\xi_i) - I(s))].
$$
We conclude using the inequality (3.9) (and further inequality (3.8)
$$
||[I_n(t)-I(t)] - [I_n(s) - I(s)]||G(\overline{\psi_g}) \le  C \ d_{\psi_g}(t,s).\eqno(4.22)
$$
Since the exact value of constant $ C $ is not essential, we get to the assertion
of theorem 4.1. \par
 The proposition of theorem 4.2 provided analogously; instead the
 $ d_{\psi}(\cdot,\cdot) $  distance we will use the metric
$ d_{\phi}(\cdot,\cdot). $ \par

\vspace{2mm}

\section{ Confidence region for solution of integral equations}
\vspace{2mm}

 Let us return to the source integral equation (1.1). We retain the  notations of the
 sections 1 and 2: $ n, \epsilon, N = N(\epsilon), \theta(m), n(m) = \theta(m) \ n,
 R(x), d(t,s), y^{(N)}_{\hat{n}}(t),y^{(N)}(t)   $ etc. \par
  Another notations:

  $$
   \psi_R(p) = |R(\cdot)|_p = \left[\int_T |R(x)|^p \ \mu(dx) \right]^{1/p}
  \eqno(5.1)
  $$
and suppose $ \psi_R(\cdot) \in \Psi,  $ i.e. $ \psi_R(b) < \infty $  for some
$ b > 2; $\par

$$
 \hat{v}_{\psi,m} = \inf_{y, y \in (0,1),\psi_R(1/y) < \infty}
\left[xy + \log \left( \frac{\psi_R^m(1/y) }{ y \cdot \log(1+1/y)} \right) \right],
\eqno(5.2)
$$

$$
\hat{v}_{\psi} = \hat{v}_{\psi,N} = \hat{v}_{\psi,N(\epsilon)},
$$

$$
\hat{I}(\epsilon) = \int_0^1 \hat{v}_{\psi,N}(H(T,d,x)) \ dx; \eqno(5.3)
$$

$$
Z_m(t,s) = \int_{T^m} K(t,x_1)K(s,x_1) K^2(x_1,x_2) K^2(x_2,x_3) \ldots
K^2(x_{m-1},x_m) \cdot
$$

$$
 f^2(x_m) \ \mu(dx_1) \mu(dx_2) \ldots \mu(dx_m), \eqno(5.4)
$$

$$
\hat{Z}(t,s) = \sum_{m=2}^N Z_m(t,s)/\theta(m). \eqno(5.5)
$$

Let $ \hat{X}(t) = \hat{X}_{\epsilon}(t) $ be a separable centered Gaussian field
with the covariation function
$ \hat{Z}(t,s): \ {\bf E} \hat{X}(t)=0, \ {\bf E} \hat{X}(t)\hat{X}(s) =
\hat{Z}(t,s). $ \par

\vspace{2mm}

{\bf A.} {\it We consider first of all the power level for integral equation.} \par

\vspace{2mm}

{\bf Theorem 5.1.a.} Assume in addition that {\it for some} $ \epsilon \in (0,1) $
 $ \hat{I}(\epsilon) < \infty. $ Then for such the value $ \epsilon $ the gaussian
random field $ \hat{X}(t) $ is $  D(\cdot,\cdot) $ continuous a.e. and

$$
\lim_{n \to \infty} {\bf P}
\left( \sqrt{n} \max_{t \in T}|y^{(N)}_{\hat{n}}(t)- y^{(N)}(t)| > u  \right)=
{\bf P} \left( \max_{t \in T} |\hat{X}(t)| > u  \right). \eqno(5.6)
$$

{\bf Theorem 5.1.b.} Assume in addition that {\it for arbitrary} $ \epsilon \in (0,1) $
 $ \hat{I}(\epsilon) < \infty. $ Then for all the values $ \epsilon $ the gaussian
random field $ \hat{X}(t) $ is $  D(\cdot,\cdot) $ continuous a.e. and the proposition
(5.6) holds.\par
 It is enough to prove only theorem 5.1.a. \par
{\bf Proof } consists on the using of theorem 4.1 to the each summand
$ S^m_{n(j)}[f](t). $ \par
 We need to prove only as before the weak compactness of the sequence  of random fields

$$
\zeta_n(t) =  \left( \sqrt{n} (y^{(N)}_{\hat{n}}(t)- y^{(N)}(t))  \right).
$$

 Indeed, let $ ||f||C(T) = 1 $  and let
 $ m = 1,2,\ldots, N; $  recall that
 $ \forall \epsilon \in (0,1) \ N = N(\epsilon) < \infty. $ We have:

$$
|K(t,s_1)K(s_1,s_2) K(s_2,s_3) \ldots K(s_{m-1},s_m) f(s_m)| \le
R(s_1)R(s_2) \ldots R(s_m),
$$

$$
|\left[K(t,s_1)- K(s,s_1)\right]\cdot K(s_1,s_2) K(s_2,s_3) \ldots K(s_{m-1},s_m) f(s_m)| \le
$$

$$
d(t,s) \cdot R(s_1)R(s_2) \ldots R(s_m).
$$
 By our condition, $ R(\cdot) \in G(\psi_R) $ or equally

 $$
 \int_T |R(s)|^p \ \mu(ds) \le \psi^p_R(p).
 $$
 Hence

 $$
 \int_{T^m} |K(t,s_1)K(s_1,s_2) K(s_2,s_3) \ldots K(s_{m-1},s_m) f(s_m)|^p
 \prod_{k=1}^m \mu(ds_k) \le
 $$

$$
\int_{T^m} \prod_{k=1}^m  R^p(s_k) \ \prod_{k=1}^m \mu(ds_k) =
\prod_{k=1}^m \int_T R^m(s_k) \mu(ds_k) \le \psi^{mp}(p)
$$
and analogously

$$
\int_{T^m} |\left[K(t,s_1)- K(s,s_1)\right]\cdot K(s_1,s_2) K(s_2,s_3) \ldots K(s_{m-1},s_m) f(s_m)|^p \prod_{k=1}^m \mu(ds_k) \le
$$

$$
d^p(t,s) \psi^{mp}(p).
$$
 Therefore, the random process $ \zeta(t) = \zeta_m(t)= $

  $$
  K(t,\xi_1)K(\xi_1,\xi_2) K(\xi_2,\xi_3) \ldots K(\xi_{m-1},\xi_m) f(\xi_m)
 $$
belongs to the space $ G(\psi^m) $ uniformly on  $ t \in T: $

$$
\sup_{t \in T} ||\zeta_m(t)||G(\psi^m) \le 1
$$
and
$$
||\zeta_m(t) - \zeta_m(s)||G(\psi^m) \le d(t,s).
$$
 The application of theorem 4.1  completes the proof of theorem 5.1.\par

\vspace{2mm}

{\bf B.}{\it Integral equations. Exponential level.} \par

\vspace{2mm}

Let us define the function $ \phi_{(N)}(\cdot)$ from the set $ \Phi: $
$$
\phi_{(N)}(p) = \left[ \frac{p^N}{\psi^N(p)} \right]^{-1}, \ p \ge 2,
\eqno(5.7)
$$
$$
\pi(\lambda) = \pi_N(\lambda) = \sup_{m = 1,2,\ldots} [m \phi_{(N)}(\lambda/\sqrt{m})].
\eqno(5.8)
$$
and suppose the finiteness of such a functions for some values $ p \ge 2. $ \par
 Denote
 $$
J = J(\epsilon) = \int_0^1 \chi_{\pi}(H(T,d_{\phi}, x)) dx.
 $$
\vspace{2mm}
{\bf Theorem 5.2.} Assume in addition that {\it for any} $ \epsilon \in (0,1) $
 $ J = J(\epsilon) < \infty. $ Then for such the value $ \epsilon $ the gaussian
random field $ \hat{X}(t) $ is $  D(\cdot,\cdot) $ continuous a.e. and the proposition
(5.6) holds.\par
{\bf Proof} is at the same as in theorem 5.1a and may be omitted. \par

\vspace{2mm}

\section{ Non-asymptotical approach}
\vspace{2mm}

We evaluate in this section the non-asymptotical probabilities for deviations

$$
\overline{{\bf P}}(u) = \sup_{n \ge 1} {\bf P^{(n)}}(u) =
\sup_{n \ge 1} {\bf P}(\sqrt{n} \sup_{t \in T} |I_n(t) - I(t)| > u) \eqno(6.1)
$$
and correspondingly

$$
\overline{{\bf Q}}(u) = \sup_{n \ge 1} {\bf Q^{(n)}}(u),
$$
where

$$
{\bf Q^{(n)}}(u) =
{\bf P}\left( \sqrt{n} \max_{t \in T}|y^{(N)}_{\hat{n}}(t)- y^{(N)}(t)| > u  \right)=
{\bf P} \left( \max_{t \in T} |\hat{Y}_n(t)| > u  \right), \eqno(6.2)
$$

$$
\hat{Y}_n(t):=\sqrt{n} (y^{(N)}_{\hat{n}}(t)- y^{(N)}(t)).  \eqno(6.3)
$$
{\bf A.} {\it Multiple integrals. Exponential level.} \par
 Assume for some function $ \phi = \phi(\lambda) \in \Phi $

$$
\sup_{t \in T} ||g(t,\eta) - I(t)||B(\phi) \stackrel{def}{=} \sigma_{\phi} < \infty.
\eqno(6.4)
$$
 The condition (6.4) is satisfied, e.g., for the natural choice of the function
 $ \phi(\lambda) = \phi_0(\lambda) $ and $ \sigma_{\phi_0} = 1. $\par
Recall that
$$
d_{\phi}(t,s) = ||[g(t,\eta) - I(t)] - [g(s,\eta) - I(s) ]||B(\phi).
$$

{\bf Theorem 6.1a.}  Suppose

$$
J(\phi):= \int_0^{\sigma_{\phi}} \chi_{\overline{\phi}}(H(T,d_{\phi},x) dx < \infty.
\eqno(6.5)
$$
Then  for the values $ u > 2\sqrt{ J(\phi) }  $

$$
\overline{{\bf P}}(\sigma_{\phi} u) \le 2
\exp \left(- \overline{\phi}^*(u - \sqrt{2 \ J(\phi) \ u} ) \right). \eqno(6.6)
$$

{\bf Proof} used the lemma 3.1. We get using the definition of the function
$ \overline{\phi}(\cdot): $

$$
\sup_{n = 1,2,\ldots} ||\sqrt{n}(I_n(t) - I(t)) ||B(\overline{\phi})  \le
\sigma_{\phi};
$$

$$
\sup_{n = 1,2,\ldots}
 ||\sqrt{n} \left[(I_n(t) - I(t)) -(I_n(s)-I(s)) \right]||B(\overline{\phi})
 \le d_{\phi}(t,s).
$$
This completes the proof of theorem 6.1. \par
\vspace{3mm}
{\bf B.} {\it Multiple integrals. Power level.} \par
{\bf Theorem 6.1b.}  Suppose for some function $ \psi = \psi(p) \in \Psi(2,b),
b > 2  $

$$
\sup_{t \in T} ||g(t,\eta) - I(t)||G(\psi) \stackrel{def}{=} \sigma_{\psi} < \infty.
\eqno(6.7)
$$
 The condition (6.7) is satisfied, e.g., for the natural choice of the function
 $ \psi(p) = \psi_0(p) $ and $ \sigma_{\psi_0} = 1. $\par
Recall that
$$
d_{\psi}(t,s) = ||[g(t,\eta) - I(t)] - [g(s,\eta) - I(s) ]||G(\psi).
$$
We define

$$
\overline{Z}(\psi):= \sigma_{\psi} + 9 \int_0^{\sigma_{\psi}}
v_{* \overline{\psi}} \left(\log(2N(T,d_{\psi},x)) \right) dx < \infty.
\eqno(6.8)
$$
If $  \overline{Z}(\psi) < \infty, $ then
for the values $ u > 2 \overline{Z}(\psi) $

$$
\overline{{\bf P}}(u) \le
\exp \left(- w_{\overline{\psi}}^*(u/\overline{Z}(\psi)) \right). \eqno(6.9)
$$

{\bf Proof} is at the same as the proof of theorem 6.1; it
used the lemma 3.2 instead lemma 4.1 and the definition of the function
$ \overline{\psi}(\cdot). $\par
 For instance, if

 $$
 \phi_0(\lambda) \sim \lambda^r, \lambda \ge 1, \ r = \const > 1,
 $$
and $ J(\phi) < \infty, $ then for $ u \ge 1 $

$$
\overline{\bf {P}}(u) \le
\exp \left(- C_1(r,J(\psi)) \ u^{\tilde{r} } \right),
$$

$$
\tilde{r} := \frac{\min(2,r)}{\min(2,r)-1}. \eqno(6.10)
$$
\vspace{3mm}
{\bf C.} {\it Integral equations. Power level.} \par
 Analogously to the theorem 6.1a may be proved the following two results.\par
{\bf Theorem 6.2.a.} Assume that {\it for some} $ \epsilon \in (0,1) $
 $ \hat{I}(\epsilon) < \infty. $  Denote

 $$
 \psi_{(N)}(p) = C_0^{-1} \frac{p \ \psi^N(p)}{\log p}, \ p \ge 2.
 $$
We deduce:

 $$
 \overline{{\bf Q}}(u) \le
 \exp \left( -w^*_{\psi_{(N)}}(C_1^{-1}u/(1+\hat{I}(\epsilon))) \right). \eqno(6.11)
 $$

\vspace{2mm}
If  for example $ X = R^d, $

$$
\mu \{x: R(x) > u \} \le \exp \left( - C_2 u^{1/\beta}  \right), u \ge 1,
$$

and  $ T  $ is  bounded open subset $ R^d,  $

$$
d_{\psi}(t,s) \le C_3 \left( \min \left(|\log|t-s|\ |^{-\gamma},1 \right)  \right),
$$

$$
\gamma = \const > \tilde{\beta} \stackrel{def}{=} \beta N(\epsilon) + 1,
$$

then $ \overline{Z}(\psi) < \infty $ and following

$$
\overline{\bf {Q}}(u)  \le
\exp \left( - C_4(C_2, C_3, \gamma,\beta,d) \ u^{1/\tilde{\beta}}  \right), u \ge 1. \eqno(6.12)
$$

\vspace{2mm}

 {\bf Remark 6.1.}  We accept that the case $ \beta = 0 $ is equivalent the boundedness
 of the function  $ R(\cdot): $
 $$
 \beta = 0 \Leftrightarrow  \vraisup_x R(x) < \infty.
 $$

\vspace{3mm}

{\bf D.} {\it Integral equations. Exponential level.} \par
We define the function $ \phi_{(N)}(\cdot)$ from the set $ \Phi: $
$$
\phi_{(N)}(p) = \left[ \frac{p^N}{\psi^N(p)} \right]^{-1}, \ p \ge 2,
\eqno(6.13)
$$
$$
\pi(\lambda) = \pi_N(\lambda) = \sup_{m = 1,2,\ldots} [m \phi_{(N)}(\lambda/\sqrt{m})].
\eqno(6.14)
$$
and suppose the finiteness of such a functions for some values $ p \ge 2. $ \par
 Denote
 $$
J = \int_0^1 \chi_{\pi}(H(T,d_{\phi}, x)) dx
 $$
and suppose $ J < \infty. $ \par
Proposition:

$$
\overline{\bf {Q}}(u)  \le \exp \left(- \pi^*(u - \sqrt{2  J u}) \right),
u > 2 J.
$$

\vspace{3mm}
 As an application: solving the equation

 $$
 \overline{\bf {P}}(u_P(\delta)) = \delta,
 $$
or correspondingly
$$
 \overline{\bf {Q}}(u_Q(\delta)) = \delta
$$
relative the variable $ u = u_P(\delta) $ or $ u = u_Q(\delta)$
where as before $ 1 - \delta, \ \delta = 0.05, \ 0.01  $ etc., is the
reliability of {\it non-asymptotic}
confidence region in the uniform norm, we conclude that with
probability at least $ 1 - \delta $

$$
 \sup_{t \in T} |I_n(t) - I(t)| \le u_P(\delta)/\sqrt{n}, \eqno(6.12)
$$

$$
  \sup_{t \in T} |y^{(N)}_{\hat{n}}(t)- y^{(N)}(t)| \le u_Q(\delta)/\sqrt{n}.
 \eqno(6.13)
$$

\vspace{2mm}

\section{ Examples}

\vspace{2mm}

 We suppose in this section that $ T $ is bounded closed domain in the space $ R^d $ with
positive Lebesgue measure $ \mu(D) = \int_D dx. $ Denote as ordinary by $ |t-s|  $
the Euclidean distance between a two points $ t,s; t,s \in T. $ \par
We assume again  $ r(S) < 1, \ r(U) < 1. $ \par
\vspace{2mm}
{\bf Example 7.1. Multiple parametric integral. "Power" level.} Recall that
$$
I(t) = \int_X g(t,x) \ \nu(dx); \ I_n(t) = n^{-1}\sum_{i=1}^n g(t,\eta(i)), \ t \in T
$$
and $ \Law(\eta(i)) = \nu.  $\par
 Assume that for some $ \alpha \in (0.1] $
 $$
 |g(t,x) - g(s,x)| \le |t-s|^{\alpha} \ Q(x) \eqno(7.1)
 $$
where for some $ \delta > 0 $

$$
\int_X Q^{d/\alpha + \delta}(x) \ \nu(dx) < \infty \Leftrightarrow
Q(\cdot) \in L_{d/\alpha+\delta}. \eqno(7.2)
$$

 Since

 $$
 N(T, |t-s|^{\alpha}, \epsilon) \asymp  \epsilon^{-d/\alpha}, \ \epsilon \to 0+
 $$
we conclude that for the value $ p_0, $ where

$$
\frac{1}{p_0} = d/\alpha + \delta,
$$
the Pizier's condition of theorem 4.1 is satisfied.\par

\vspace{2mm}
{\bf Example 7.2. Multiple parametric integral. "Exponential" level.} Here we suppose

$$
 |g(t,x) - g(s,x)| \le
 \max[|\log|t-s||^{-\gamma},1] \ Q(x), \eqno(7.3)
 $$
where

$$
\nu \{x: Q(x) > u \} \le  C_1 \exp \left( -C_2 u^{1/\beta}  \right), \eqno(7.4)
$$
and
$$
\beta, \gamma = \const, \ \gamma > \beta.  \eqno(7.5)
$$
We conclude that the condition of theorem 4.1 is satisfied.\par
 Note that the condition (7.3) is equivalent to the following inequality:

 $$
 \sup_{p \ge 1} |Q|_p/p^{\beta} < \infty,
 $$
or equally

$$
Q(\cdot) \in G \left(\psi_{\beta} \right), \ \psi_{\beta}(p) \stackrel{def}{=}
p^{\beta}.
$$

The condition of theorem 4.1 is satisfied also when

 $$
 \sup_{p \ge 1} |Q|_p/ \left[p^{\beta} \ (\log p)^{\beta_2} \right] < \infty, \
 \beta_2 = \const > 0.
 $$

\vspace{2mm}
{\bf Example 7.3. Integral equation. "Power" level.} \par
{\bf A.} Theorem 5.1a. \par
Let again  $ T $ is bounded open subset of the space $ R^d $ and assume as before

$$
\mu \{x: R(x) > u \} \le \exp \left( - C_1 u^{1/\beta}  \right), u \ge 1,
\eqno(7.6)
$$

$$
d_{\psi}(t,s) \le C_2 \left( \min \left(|\log|t-s|\ |^{-\gamma},1 \right)  \right),
\eqno(7.7)
$$

$$
\gamma = \const >  \beta N(\epsilon) + 1. \eqno(7.8)
$$
 Then all the conditions of theorem 5.1a are satisfied.
 In particular, $ \overline{Z}(\psi,\epsilon) < \infty. $\par

\vspace{3mm}

{\bf B.} Theorem 5.1b. \par
Let the condition (7.6) be satisfied.  Suppose also  (instead conditions (7.7) and
(7.8) )

$$
d(t,s) \le C_4 |t-s|^{\alpha}, \ \alpha = \const \in (0,1].
$$
Then all the conditions of theorem 5.1b are satisfied.\par
In particular, $ \forall \epsilon \in (0,1) \
\overline{Z}(\psi,\epsilon) < \infty. $\par

\vspace{2mm}

\section{ Derivative computation}

\vspace{2mm}

 Let us return to the source equation (1.1).  We consider in this section the case
 $ T = [0,1] $ with the classical Lebesgue measure $ \mu(dx) = dx. $ Suppose
 $ f $ be continuous differentiable: $ f^/ \in C[0,1] $  and that there exists  a
 continuous relative the variable $ t $ function
$$
V(t,s) = \frac{\partial K(t,s)}{\partial t}.
$$
 We obtain denoting $ Y(t) = y^/(t) $ after the differentiation of equation (1.1)

$$
Y(t) = f^/(t) + \int_0^1 V(t,s) y(s) ds \stackrel{def}{=} \tilde{f}(t) + V[y](t). \eqno(8.1)
$$
{\it  We do not entrust the contraction condition on the kernel } $ V = V(t,s). $\par
 The solution $ Y(t) $ may be written as follows.
 $$
 Y(t) = \tilde{f}(t) + \sum_{m=1}^{\infty}\sigma_m(t),
 $$
 where

$$
\sigma_1(t) = \int_0^1 V(t,x) f(x) dx,  \ m = 2,3,\ldots \Rightarrow \sigma_m(t)=
$$

$$
\int_0^1ds \int_{T^m} V(t,s) K(s,s_1) K(s_1,s_2) \ldots K(s_{m-1},s_m) f(s_m) \ ds
\ ds_1 ds_2 \ldots ds_m. \eqno(8.2)
$$

Let the number $ \epsilon, \ \epsilon \in (0,1) $ be a given. We retain here the
notations of the second section: $ N = N(\epsilon), \ n(m) = \theta(m) \cdot N(\epsilon) $ etc. \par
 We choose as the deterministic approximation $ Y^{(N)}(t) $ for $ Y(t) $ as before
 the expression

$$
 Y^{(N)}(t) = \tilde{f}(t) + \sum_{m=1}^{N(\epsilon)}\sigma_m(t). \eqno(8.3)
$$

The accuracy of the approximation (8.3) ("bias") in the uniform norm
may be estimated as follows.

$$
\sup_{t \in [0,1]} |Y^{(N)}(t) - Y(t)| \le ||V|| \cdot
\sup_{t \in [0,1]} |y^{(N)}(t) - y(t)| \le ||V|| \cdot \epsilon. \eqno(8.4)
$$
 Recall that

 $$
 ||V|| = \sup_{t \in [0,1]} \int_0^1 |V(t,x)| dx < \infty.
 $$
 Further, we offer for the $ Y^{(N)}(t) $ calculation the following Monte-Carlo approximation.  Let us consider separately the expression for $ \sigma_m(t), \
 m = 2,3,\ldots, N. $  Let the random variable  $\zeta $ be uniform distributed
 in the set $ [0,1] $ and let the random vector
 $ \vec{\xi}= \{\xi_1, \xi_2, \ldots, \xi_m \} $
 be uniform distributed in the multidimensional  cube $ [0,1]^m; $ let also
 $ \{\zeta^{(j)} \}, \ \{\vec{\xi^{(j)} } \}= \{\xi_1^j, \xi_2^j, \ldots, \xi_m^j \},
 \ j = 1,2,\ldots, n(m) $ be independent copies of the vector $ \{\zeta, \vec{\xi}\}.$ \par
 Note that

 $$
\sigma_m(t) = {\bf E} V(t,\zeta) K(\zeta,\xi_1) K(\xi_1,\xi_2) \ldots
K(\xi_{m-1}, \xi_m) f(\xi_m). \eqno(8.5)
$$
 We can offer therefore the following Monte-Carlo approximation $ \hat{\sigma}_m(t) $
 for $ \sigma_m(t): $

$$
\hat{\sigma}_m(t) := \frac{1}{n(m)} \sum_{j=1}^{n(m)}V(t,\zeta^{(j)}) f(\xi^{(j)}_m) \times
$$

$$
K(\zeta^{(j)},\xi^{(j)}_1) K(\xi^{(j)}_1, \xi^{(j)}_2) \ldots K(\xi^{(j)}_{m-1}, \xi^{(j)}_m) \eqno(8.6)
$$
and correspondingly the following approximation $ Y^{(N)}_n(t) $
for $ Y^{(N)}(t): $

$$
 Y^{(N)}_n(t) = \tilde{f}(t) + \sum_{m=1}^{N(\epsilon)} \hat{\sigma}_m(t). \eqno(8.7)
$$
 Here $ n $ denotes as in the section 2 the common quantity  of elapsed random variables.\par
 Note that

 $$
 {\bf Var}_m \stackrel{def}{=} {\bf Var}[\hat{\sigma}_m(t)] \le \frac{1}{n(m)}
\int_0^1 ds \int_{T^m}  ds_1 ds_2 \ldots ds_m
 V^2(t,s) f^2(s_m) \times
 $$

 $$
  K^2(s,s_1) K^2(s_1,s_2) \ldots K^2(s_{m-1},s_m) \le \frac{1}{n(m)} ||V^{(2)}|| \
 \cdot \ ||U^m||.  \eqno(8.8)
 $$
 It is easy to calculate as in the second section that

$$
\sup_{t \in T}{\bf Var}[Y^{(N)}_n(t)] =  
\sum_{m=1}^{N(\epsilon)} \sup_{t \in T} {\bf Var}_m \le C/n. \eqno(8.9)
$$
 The inequality (8.9) show us that the speed of convergence $ Y^{(N)}_n(t) $ to
$ Y^{(N)}(t) $ is equal to $ 1/\sqrt{n} $ in each fixed point $ t_0 \in T. $ \par
 In order to establish this result in the uniform norm for the sequence of random
 fields $ Y^{(N)}_n(t), $  we need to introduce some new notations and conditions. \par

$$
R_V(x) = \vraisup_{t \in T} |V(t,x)|, \ d_V(t,s) =
\vraisup_{x \in T} \frac{|V(t,x)-V(s,x)|}{R_V(x)},
$$

$$
\psi_V(p) = |R_V(\cdot)|_p,
$$
so that

$$
|V(t,x)| \le R_V(x), \ |V(t,x)-V(s,x)| \le d_V(t,s) \ R_V(x), \eqno(8.10)
$$

$$
\exists b_V > 2, \ \forall p < b_V \  \psi_V(p) < \infty;
$$
 $$
 b_{RV} = \min(b,b_V) = \const > 2;
 $$

$$
V_1(t,s)= \int_T V(t,z) V(s,z) dz;
$$

$$
V_m(t,s)=\int_0^1 dz \int_{T^m} V(t,z)V(s,z) K^2(z,x_1) K^2(x_1,x_2) \ldots
K^2(x_{m-1}, x_m) \times
$$

$$
f^2(x_m) \ dx_1 dx_2 \ldots dx_m, \ m = 2,3,\ldots,N(\epsilon); \eqno(8.11)
$$

$$
Z_V(t,s) = Z_{V,N}(t,s)= \sum_{m=1}^{N(\epsilon)} \frac{1}{\theta(m)} V_m(t,s);
$$

$$
\Theta_V(p) = \psi_V(p) \cdot \left[\psi_R(p) \right]^N, \ p \in (2,b_{RV});
$$

$$
Z_V:= 1 + 9\int_0^1 v_{* \Theta_V} (H(T, d_V,x)) \ dx; \eqno(8.12)
$$

 Let $ X_V(t) = X_{V, \epsilon}(t) $ be a separable centered Gaussian field
with the covariation function
$ Z_V(t,s): \ {\bf E} X_V(t)=0, \ {\bf E} X_V(t)X_V(s) = Z_V(t,s). $\par

\vspace{4mm}

{\bf Theorem 8.1.} (Power level.)
 Let $ Z_V < \infty. $  Then the limiting Gaussian random field
 $ X_V(\cdot) $ is $ d_{V} $ continuous a.e. and

\vspace{2mm}

$$
\lim_{n \to \infty} {\bf P}(\sqrt{n} \sup_{t \in T}|Y^{(N)}_n(t)-Y^{(N)}(t)| > u) =
{\bf P}(\sup_{t \in T} |X_V(t)| > u), \ u > 0; \eqno(8.13)
$$

$$
\sup_n {\bf P}(\sqrt{n} \sup_{t \in T}|Y^{(N)}_n(t)-Y^{(N)}(t)| > u) \le
\exp \left(- w^*_{\Theta_V}(u/Z_V) \right), \ u > 2 Z_V. \eqno(8.14)
$$
{\bf Proof } is at the same as the proof of theorems 5.1a and 5.1b.\par
Namely, let us introduce the following random processes  (fields)

$$
\zeta_{m,V}(t) = V(t,\zeta) K(\zeta,\xi_1) K(\xi_1,\xi_2)\ldots K(\xi_{m-1},\xi_m) f(\xi_m)-
$$

$$
V\cdot S^m[f](t), \ m=1,2,\ldots, N \eqno(8.15)
$$
so that $ {\bf E}\zeta_{m,V}(t)=0, $ and its independent copies

$$
\zeta_{m,V}^{(j)}(t) = V(t,\zeta^{j}) K(\zeta^{(j)},\xi_1^{(j)}) K(\xi_1^{(j)},\xi_2^{(j)})\ldots K(\xi_{m-1}^{(j)},\xi_m^{(j)}) f(\xi_m^{(j)})-
$$

$$
V\cdot S^m[f](t), \ m=1,2,\ldots, N; \eqno(8.16)
$$

$$
\Xi_m^{(n)}(t) = \frac{1}{\sqrt{n(m)}} \sum_{j=1}^{n(m)} \zeta_{m,V}^{(j)}(t);
\eqno(8.17)
$$

 $$
 \zeta_{V}^{(j)}(t) = \zeta_{N,V}^{(j)}(t), \ \Xi^{(n)}(t) = \Xi_N^{(n)}(t).
 $$
 It is sufficient to consider only the case $ m = N. $ We need to prove the tightness
 the random fields $  \Xi^{(n)}(t). $ \par
  We have assuming without loss of generality $ \sup_t |\tilde{f}(t)| = 1: $

 $$
|\zeta_{m,V}(t)| \le R_V(\zeta) \ R(\xi_1) R(\xi_2) \ldots R(\xi_m),
 $$

 $$
|\zeta_{m,V}(\cdot)|_p \le \psi_V(p) \ \psi_R^m(p),
 $$
and using the Rosenthal's inequality

$$
|\Xi^{(n)}(t)|_p \le C_0^{-1} p \ \psi_V(p) \ \psi_R^m(p) /\log p. \eqno(8.18)
$$
Analogously

$$
|\Xi^{(n)}(t) - \Xi^{(n)}(s)|_p \le d_V(t,s) \cdot
C_0^{-1} p \ \psi_V(p) \ \psi_R^m(p) /\log p, \eqno(8.19)
$$
 We conclude after summing over $ m: $

$$
\sup_n \sup_t || \left[\sqrt{n}(Y^{(N)}_n(t)-Y^{(N)}(t)) \right]|| G(\Theta_V) \le 1, \eqno(8.20)
$$

$$
\sup_n
|| \left[\sqrt{n}(Y^{(N)}_n(t)-Y^{(N)}(t)) \right] -
\left[\sqrt{n}(Y^{(N)}_n(s)-Y^{(N)}(s)) \right]|| G(\Theta_V) \le
$$
$$
 d_V(t,s). \eqno(8.21)
$$
 The assertion of theorem 8.1. follows now from lemma 3.2 and theorem 3.2.\par

\vspace{3mm}

{\bf Example 8.1.} If

$$
\mu \{x: R(x) > u \} \le \exp \left(-C u^{1/\beta}  \right), \ \beta = \const \ge 0,
$$

$$
\mu \{x: R_V(x) > u \} \le \exp \left(-C u^{1/\omega}  \right), \ \omega = \const \ge 0,
$$

$$
d_V(t,s) \le C \left[\min \left( |\log|t-s| \ |^{-\gamma}  \right),1 \right], \gamma =
\const,
$$
and

$$
\gamma > N(\epsilon) \beta + \omega + 1,
$$
 then the conditions of theorem 8.1 are satisfied. As a consequence:

$$
\sup_n {\bf P}(\sqrt{n} \sup_{t \in T}|Y^{(N)}_n(t)-Y^{(N)}(t)| > u) \le
\exp \left(- C(\beta,\gamma)u^{1/(N \beta + \omega + 1)}  \right), \ u > 2.
\eqno(8.22)
$$

In the case when
$$
d_V(t,s) \le C |t-s|^{\alpha}, \ \alpha = \const > 0,
$$
then the conditions of theorem 8.1 are satisfied for arbitrary values $ \epsilon
\in (0,1). $ \par

\vspace{3mm}
{\bf Remark 8.1.} CLT in the space $ C^1[T]. $\par
\vspace{3mm}

 We used in this section in fact the Central Limit Theorem in the space of continuous
differentiable function $ C^1[T]. $ \par

\vspace{3mm}
{\bf Remark 8.2.} Some generalizations.
\vspace{3mm}
 Let $ A = A_t $ be any linear operator, not necessary to be bounded, defined on some
 (dense or not) subspace of the space $ C(T), $  for example, differential operator
 $ A = d/dt, $ partial differential operator, Laplace's operator  etc. \par
  We write formally

 $$
 A_t y = A_t f(t) + \int_T A_t K(t,s) y(s) \mu(ds).
 $$
Using at the same considerations, we might obtain the CLT in the space $ C_A[T] $
consisting on the continuous functions $ g = g(t) $ with  continuous $ Ag(t) $
equipped by the "energy" norm

$$
||g||_A = \sup_{t \in T} |g(t)| + \sup_{t \in T} |Ag(t)|.
$$
 We can conclude as a consequence that the rate of convergence of the random approximation for the Monte-Carlo approximation

$$
 Y_A^{(N)}(t) := Af(t) + \sum_{m=1}^{N(\epsilon)}A\sigma_m(t)
$$
in the space $ C_A[T] $ is equal to $ 1/\sqrt{n} $
and the bias is less than $ C\cdot \epsilon: $ \par

$$
|| Ay - Y_A^{(N)}(t)||_A  \le  \epsilon \cdot \sup_{t \in T} \int_T |A_t K(t,s)| \ ds.
\eqno(8.23)
$$

\vspace{3mm}

We consider now the exponential level for the derivative computation.  Recall
(see section 3) that the estimations through $ B(\phi) $ spaces (exponential level)
have advantage in comparison  to the estimations using the $ G(\psi) $ technique
if the correspondent $ \phi(\cdot) - $ function there exists.\par
 Define a function
$$
\phi_V(p) = \left[\frac{p}{\Theta_V(p)}\right]^{-1}, \eqno(8.24)
$$
if there exists for some interval $ p \in [1,b], \ b = \const \ge 2, $ and suppose
$ \phi_V(\cdot) \in \Phi, $ and introduce the correspondent distance

$$
d_{\phi_V}(t,s) = ||\zeta_{N,V}(t) - \zeta_{N,V}(s)||B(\phi_V). \eqno(8.25)
$$

 Moreover, assume that the following integral converges:

 $$
 K_V:= \int_0^1 \chi_{\overline{\phi_V}} \left( H \left(T, d_{\phi_V}, x \right)  \right)
  \ dx < \infty. \eqno(8.26)
 $$

\vspace{2mm}

{\bf Theorem 8.2.} (Exponential level). \par
\vspace{2mm}

 Let $ K_V < \infty. $
 Then the limiting Gaussian random field
 $ X_V(\cdot) $ is $ d_{\phi_V} $ continuous a.e. and

\vspace{2mm}

$$
\lim_{n \to \infty} {\bf P}(\sqrt{n} \sup_{t \in T}|Y^{(N)}_n(t)-Y^{(N)}(t)| > u) =
{\bf P}(\sup_{t \in T} |X_V(t)| > u), \ u > 0; \eqno(8.27)
$$

$$
\sup_n {\bf P}(\sqrt{n} \sup_{t \in T}|Y^{(N)}_n(t)-Y^{(N)}(t)| > u) \le
\exp \left(- \phi_{V}(u - \sqrt{2 K_V u} ) \right), \ u > 2 K_V. \eqno(8.28)
$$
{\bf Proof } is at the same as the proof of theorems 8.1 and may be omitted\par

\vspace{2mm}

{\bf Remark 8.3} As we know in the remark 4.1, the tail functions for the maximum
distributions  in (5.6) and (8.27), i.e. probabilities

$$
{\bf P}(\max_{t \in T} |\hat{X}(t)| > u), \ {\bf P}(\max_{t \in T} |X_V(t)| > u)
$$
have the exact asymptotic of a view, e.g.,

$$
{\bf P}(\max_{t \in T} |\hat{X}(t)| > u) \sim  C(X,T) \ u^{\kappa - 1} \exp \left(-u^2/(2\sigma_+^2) \right),
$$
where $ \sigma_+^2 = \max_t {\bf Var} [\hat{X}(t)] $ etc.\par
 The consistent estimations of the parameters $ \kappa, \sigma^2_+ $ and
 $ C(\hat{X},T) $ may be the implemented as in the remark 4.2, where we can computed
 multiple integrals by Monte-Carlo approximation.\par
 Obviously, these parameters may be estimated through the analytical  expression for
 the functions $ f,g, K. $ \par
\vspace{2mm}
{\bf Example 8.2.} The conditions of theorem 8.2 are satisfied if for example both
the r.v. $ R(x), \ R_V(x) $ are essentially bounded,  $ X $ is bounded subset
of the space $ R^d $ and

$$
d_{\phi_V}(t,s) \le C_5 \left( \min \left(|\log|t-s|\ |^{-\gamma},1 \right)  \right),
\ \gamma = \const > 2 \eqno(8.29)
$$
or moreover in the case when (instead the condition (8.29))
$$
d_{\phi_V} (t,s) \le C_6 \ |t-s|^{\alpha}, \ \exists \alpha = \const \in (0,1].
\eqno(8.30)
$$

\vspace{2mm}

\section{ Concluding remarks}
\vspace{2mm}

{\bf A. Another method.} \par
Let us consider the following equivalent modification of source equation (1.1):

$$
 y_{\lambda}(t) = f(t) +  \lambda \int_T K(t,s) \ y(s) \ \mu(ds) = f(t) +
 \lambda S[y](t), \eqno(9.1)
$$
where $ \lambda = \const \in (0,1) $ and as before

$$
\rho_1 = r(S) < 1, \ \rho = r(U) = r \left(S^{(2)} \right) < 1
$$
and suppose that the function $ f(\cdot) $ and the kernel $ K(\cdot, \cdot) $ satisfies
all the conditions of the section 1. \par

 The solution $ y_{\lambda} = y_{\lambda}(\cdot) $ may be written as follows:

 $$
 y_{\lambda} = f + \sum_{m=1}^{\infty} \lambda^m S^m[f]. \eqno(9.2)
 $$

Let us introduce the so-called {\it geometrical distributed} integer random  variables
$ \tau: $

$$
{\bf P} (\tau = m) = (1-\lambda) \lambda^m, \ m = 0,1,2,\ldots;
$$
then

$$
(1-\lambda) y_{\lambda} = {\bf E} S^{\tau}[f].  \eqno(9.3)
$$
 Let $ M = 2,3,\ldots $ be  arbitrary integer number and let
$ \ \tau_j, \ j=1,2,\ldots,N $ be independent copies of $ \tau. $ It may be offered
as an consistent approximation for $ (1-\lambda) y_{\lambda}$ the following expression:

$$
(1-\lambda) y_{\lambda,n}(t)= M^{-1} \sum_{j=1}^M S_n^{\tau(j)},
$$
where the value $ S_n^{\tau(j)} $ has a following Monte-Carlo approximation:

$$
S_n^{\tau(j)} \approx \hat{S}_n^{\tau(j)} \stackrel{def}{=}(n(j))^{-1}
\sum_{i=1}^{n(j)} \vec{K}^{(\tau(j))}[f](\vec{\xi}^i_{t,\tau(j)}).
$$
So,

$$
(1-\lambda) y_{\lambda,n}= M^{-1} \sum_{j=1}^M S_n^{\tau(j)} (n(j))^{-1}
\sum_{i=1}^{n(j)} \vec{K}^{(\tau(j))}[f](\vec{\xi}^i_{t,\tau(j)}).\eqno(9.4)
$$

We denote also

$$
\tilde{n} := (n(1), n(2), \ldots, n(M));
$$
$ \tilde{n} $  is any $ M - $ tuple of integer positive numbers.\par
 We obtain after some calculations:

 $$
\phi(M,\tilde{n}) :=
 {\bf Var } ( y_{\lambda,n}) \asymp \frac{1}{M}  \sum_{j=1}^M  \frac{1}{n(j)}.
 $$

Note as before that the common  amount of used $ T - $ valued random variables is
equal to

$$
A(M,\tilde{n}):= M \cdot \sum_{j=1}^M  j \ n(j).
$$
We conclude solving the following constrained extremal problem:

 $$
\phi(M,\tilde{n}) \to \min  /A(M,\tilde{n}) = n,
 $$
that the minimal value of the function $ \phi(M,\tilde{n}) $  under condition
$ A(M,\tilde{n}) = n $ is asymptotically as $ n \to \infty $ equivalent to

$$
\min \phi(M,\tilde{n}) \  /[A(M,\tilde{n}) = n] \ \asymp n^{-1/2}. \eqno(9.5)
$$
Thus, the optimal speed of convergence $ y_{\lambda,n}  $ to the solution $ y $ is
asymptotical equal to $ n^{-1/4}, $ in contradiction to the first offered method.\par

\vspace{3mm}

 Another {\it sufficient} conditions for CLT in the space of continuous functions
$ C(T) $ see, e.g. in the works \cite{Gine1}, \cite{Jain1}, \cite{Pizier2}, \cite{Pollard1}. \par
 We adapt here only the results belonging to Marcus M.B. and Jain N.C. for the integral
 $ I(t) = \int_X g(t,x) \nu(dx) $ calculation. Namely, if

 $$
 |g(t,x) - g(s,x)| \le M(x) \cdot \rho(t,s), \eqno(9.6)
 $$
(condition of factorization),
$$
\int_X M^2(x) \nu(dx) < \infty, \eqno(9.7)
$$
(moment condition),
where $ \rho(t,s) $ is some metric on the set $ T $  for which

$$
\int_0^1  H^{1/2}(T,\rho, z) \ dz < \infty, \eqno(9.8)
$$
(entropy condition),
then the sequence of random processes (fields) $ \sqrt{n} (I_n(t) - I(t)) $  satisfies
the Central Limit Theorem in the space $ C(T,\rho). $ \par
 Analogously, if

 $$
 \int_X \exp \lambda \left( [g(t,x) - g(s,x)] - [ I(t) - I(s) ]  \right) \nu(dx) \le
 $$
 $$
 \exp(0.5 A^2 \lambda^2 \tau^2(t,s)), \ A = \const > 0, \eqno(9.9)
 $$
 or equally

 $$
 || \left( [g(t,x) - g(s,x)] - [ I(t) - I(s) ]  \right) ||\sub \le A \tau(t,s)
 $$
  (subgaussian condition), for some (semi-) distance $ \tau(\cdot,\cdot)$ for which

$$
\int_0^1 H^{1/2}(T,\tau, x) dx < \infty, \eqno(9.10)
$$
(entropy condition), then also  the sequence of random processes (fields) $ \sqrt{n} (I_n(t) - I(t)) $  satisfies the Central Limit Theorem in the space $ C(T,\tau). $ \par

\vspace{2mm}

\section{ Necessity of CLT conditions }
\vspace{2mm}

We will confer in this section the necessity of some conditions for the Central
Limit Theorem in the space of continuous functions $ C(T,\rho). $ \par
{\bf A. Condition of factorization.}\par
Let $ \zeta(t) $  be continuous with probability one random field relative some
distance. Then there exist a non-random continuous distance $ \rho = \rho(t,s) $
and a random variable $ M = M(\omega) $ for which

$$
|\zeta(t) - \zeta(s)| \le M(\omega) \cdot \rho(t,s),  \eqno(10.1)
$$
see \cite{Ostrovsky19}, \cite{Buldygin3}. \par
{\bf B. Moment condition.} \par
 Assume in addition the field $ \zeta(t) $ satisfies the Central Limit Theorem in
 the space of continuous functions. Then $ \zeta(\cdot) $ has the weak second moment. Following, see \cite{Ostrovsky19},

 $$
 {\bf E} M^2 < \infty. \eqno(10.2)
 $$

{\bf C. Entropy condition.}\par
 Let the conditions (9.6) and (9.7) for the field $ \zeta = \zeta(t) $
  be satisfied; assume also the field $ \zeta(t) $
satisfies the Central Limit Theorem in the space of continuous functions. \par
 Suppose in addition $ T = [0,2\pi]^d, \ d = 1,2,\ldots $ and that the distance
$ \rho(t,s) $ from the inequality(10.1) is translation invariant:

$$
\rho(t,s) = \rho(|t-s|), \eqno(10.3)
$$
where $ t \pm s = t \pm s (\mod 2 \pi).$ \par
 Other assumption. Let the random field $ \zeta(t) $
satisfies the Central Limit Theorem in the space of continuous functions. \par
 Denote by $ \overline{\zeta}_{\infty}(t) $ the separable centered Gaussian field
with at the same covariation function as $ \eta(t): $

$$
{\bf E}\overline{\zeta}_{\infty}(t)\overline{\zeta}_{\infty}(s) =
{\bf E} \zeta(t) \zeta(s), \eqno(10.4)
$$
then

$$
\rho_{\infty}(t,s):=
{\bf E}^{1/2}\left[\overline{\zeta}_{\infty}(t)- \overline{\zeta}_{\infty}(s)\right]^2 =
$$

$$
{\bf E}^{1/2}\left[\zeta(t)- \zeta(s)\right]^2 =: d_{\zeta}(t,s). \eqno(10.5)
$$
 It follows from the condition  (10.2) that

 $$
 d_{\zeta}(t,s) \le C \rho(|t-s|);
 $$
assume in addition the conversely:

 $$
\rho(|t-s|) \le C_1 d_{\zeta}(t,s). \eqno(10.6)
 $$
As long as the distance $ \rho_{\infty}(t,s) $ is linear equivalent to the distance
$ \rho(|t-s|): $

 $$
 d_{\zeta}(t,s) \asymp \rho(|t-s|), t,s \in T,
 $$
the convergence of entropy integral (9.8) follows immediately from the
famous result of X.Fernique \cite{Fernique1}. \par
{\bf D. Subgaussian condition.}\par
 Let the condition (9.9) be satisfied. Suppose alike in the last pilcrow
in addition $ T = [0,2\pi]^d, \ d = 1,2,\ldots $ and that the distance
$ \tau(t,s) $ from the inequality (9.9) is translation invariant:

$$
\tau(t,s) = \tau(|t-s|),
$$
where $ t \pm s = t \pm s (\mod 2 \pi).$ \par
 Then the (weak) limiting Gaussian field $ \nu = \nu(t) $ in the CLT for the space of continuous functions satisfies at the same condition (9.9). \par
  It is evident that

$$
|\nu(t) - \nu(s)|_2 \le \tau(|t-s|);
$$
assume conversely, i.e. that

$$
\tau(|t-s|) \le C |\nu(t) - \nu(s)|_2;
$$
then the convergence of entropy integral (9.10) follows again from the
result of X.Fernique \cite{Fernique1}. \par

\end{document}